\documentclass[11pt,leqno]{article}
\usepackage{amssymb}
\usepackage{amsmath}
\usepackage{latexsym}
\usepackage{amsfonts}
\usepackage{epsfig}

\begin{document}
\newcommand{\st}{{\,|\,}}		
\renewcommand{\frak}{\mathfrak}		
\renewcommand{\Bbb}{\mathbb}		
\newcommand{\ds}{\displaystyle}		
\newcommand{\ul}{\underline}		
\newcommand{\simga}{\sigma}		
\newcommand{\lamdba}{\lambda}		
\newcommand{\tua}{\tau}			
\newcommand{\tensor}{\otimes}		
\newcommand{\surjection}{\twoheadrightarrow}
\newcommand{\union}{\cup}		
\newcommand{\intersect}{\cap}		
\newcommand{\imbedding}{\hookrightarrow}
\newcommand{\gdn}{G_{d,n}}		
\newcommand{\idn}{I(d,n)}		
\newcommand{\svw}{S^v_w}		
\newcommand{\svwm}{\svw(m)}		
\newcommand{\sv}{S^v}			
\newcommand{\ui}{{\underline{i}}}	
\newcommand{\pos}{{\frak N}}		
\newcommand{\smvw}{\begin{mbox}
		{\rm SM}\end{mbox}^v_w}	
\newcommand{\smvwm}{\smvw(m)}		
\newcommand{\mon}{\frak S}		
\renewcommand{\part}{\frak B}		
\newcommand{\vdeg}[1]{\mbox{{\rm $v$-degree}}(#1)}
\newcommand{\vdegree}{\vdeg}		
\newcommand{\degree}[1]{\mbox{{\rm degree}}(#1)}
\newcommand{\product}{\prod}		
\newcommand{\tv}{T^v}			
\newcommand{\uv}{U^v}			
\newcommand{\uvm}{U^v(m)}		
\newcommand{\tvwm}{T^v_w(m)}		
\newcommand{\smvtopv}{\begin{mbox}
		{\rm SM}\end{mbox}^{v,v}} 
\newcommand{\smvtopvm}{\begin{mbox}
		{\rm SM}\end{mbox}^{v,v}(m)} 
\newcommand{\smvtopvw}{\begin{mbox}
		{\rm SM}\end{mbox}^{v,v}_w} 
\newcommand{\smvtopvwm}{\begin{mbox}
		{\rm SM}\end{mbox}^{v,v}_w(m)} 
\newcommand{\smvvm}{\begin{mbox}
		{\rm SM}\end{mbox}^{v}_v(m)} 
\newcommand{\pitilde}{\tilde{\pi}}	
\newcommand{\phitilde}{\tilde{\phi}}	
\newcommand{\roots}{\frak R }		
\renewcommand{\min}{\frak T}		
\newcommand{\minwj}{{\min}^w_j}		
\newcommand{\monj}{\mon_j}		
\newcommand{\monw}{{\mon}_w}		
\newcommand{\monwj}{{\mon}_w^j}		
\newcommand{\monwk}{{\mon}_w^k}		
\newcommand{\monwone}{{\mon}_w^1}		
\newcommand{\monwtwo}{{\mon}_w^2}		
\newcommand{\piece}{\frak P}			
\newcommand{\blockc}{\blc}			
\newcommand{\block}{\frak B}		
\newcommand{\blc}{\frak C}		
\newcommand{\blcnot}{\blc_0}		
\newcommand{\blcone}{\blc_1}		
\newcommand{\blocka}{\frak A}		
\newcommand{\blockd}{\frak D}		
\newcommand{\blocke}{\frak E}		
\newcommand{\card}[1]{|#1|}		
\newcommand{\init}[1]{\mbox{{\rm in}}(#1)}
\newcommand{\poone}{>_1}		
\newcommand{\potwo}{>_2}		
\newcommand{\pothree}{>_3}		
\newcommand{\pofour}{>_4}		
\newcommand{\uvw}{{\frak U}^v_w}	
\newcommand{\nvw}{N^v_w}		
\newcommand{\rvw}{R^v_w}		
\newcommand{\path}{\Lambda}		
\newcommand{\bstart}{\beta^{\mbox{\begin{rm}start\end{rm}}}}
\newcommand{\bstartj}{\beta^{\mbox{\begin{rm}start\end{rm}}}_j}
\newcommand{\bfinish}{\beta^{\mbox{\begin{rm}finish\end{rm}}}}
\newcommand{\bfinishj}{\beta^{\mbox{\begin{rm}finish\end{rm}}}_j}
\addtolength{\itemsep}{-100mm}
\title{Hilbert functions of points on\\
	 Schubert varieties in the Grassmannian}
\author{V.~Kodiyalam  and K.N.~Raghavan \\
Institute of Mathematical Sciences \\
C.I.T.~Campus, Taramani\\
Chennai 600 113 INDIA \\
E-mail: vijay@imsc.ernet.in and knr@imsc.ernet.in
}

\maketitle


\newtheorem{theorem}{Theorem}[section]
\newcommand{\bthm}{\begin{theorem}}
\newcommand{\ethm}{\end{theorem}}

\newtheorem{proposition}[theorem]{Proposition}
\newcommand{\bpr}{\begin{proposition}}
\newcommand{\bprop}{\begin{proposition}}
\newcommand{\epr}{\end{proposition}}
\newcommand{\eprop}{\end{proposition}}

\newtheorem{definition}[theorem]{Definition}
\newcommand{\bdefn}{\begin{definition}\begin{rm}}
\newcommand{\edefn}{\end{rm}\end{definition}}

\newtheorem{example}[theorem]{Example}
\newcommand{\bexample}{\begin{example}\begin{rm}}
\newcommand{\eexample}{\ $\Box$\end{rm}\end{example}}

\newtheorem{remark}[theorem]{Remark}
\newcommand{\bremark}{\begin{remark}\begin{rm}}
\newcommand{\eremark}{\end{rm}\end{remark}}

\newtheorem{corollary}[theorem]{Corollary}
\newcommand{\bcor}{\begin{corollary}}
\newcommand{\ecor}{\end{corollary}}

\newtheorem{exercise}[theorem]{Exercise}
\newcommand{\bex}{\begin{exercise}\begin{rm}}
\newcommand{\eex}{\end{rm}\end{exercise}}

\newtheorem{lemma}[theorem]{Lemma}
\newcommand{\blem}{\begin{lemma}}
\newcommand{\blemma}{\begin{lemma}}
\newcommand{\elemma}{\end{lemma}}
\newcommand{\elem}{\end{lemma}}

\newenvironment{proof}[1][]{{\sc Proof #1}:\ }{\ $\Box$}
\newcommand{\bpf}{\begin{proof}}
\newcommand{\bproof}{\begin{proof}}
\newcommand{\epf}{\end{proof}}
\newcommand{\eproof}{\end{proof}}

\newenvironment{prooftwo}[1][]{{\sc Proof #1}:\ }{}
\newcommand{\bpftwo}{\begin{prooftwo}}
\newcommand{\bprooftwo}{\begin{prooftwo}}
\newcommand{\epftwo}{\end{prooftwo}}
\newcommand{\eprooftwo}{\end{prooftwo}}

\newenvironment{solution}{{\sc Solution}:\ }{\ $\Box$}
\newcommand{\bsol}{\begin{solution}}
\newcommand{\esol}{\end{solution}}

\section{Introduction}\label{sintro}
This paper arose out of our reading of the recent paper \cite{kl}
of Kreiman and Lakshmibai.   They consider the following problem.
Let $X$ be a Schubert variety in the Grassmannian $\gdn$ of $d$-planes
in $n$-dimensional linear space.    Given a point $x$ on $X$,  what is
the Hilbert function of the tangent cone to $X$ at $x$?

In case the point $x$ is the ``identity coset'',  Kreiman and
Lakshmibai deduce, from well known results of Hodge~\cite[Chapter XIV]{hodge} and
Musili~\cite{smt} about homogeneous co-ordinate rings of Schubert varieties in
the Grassmannian,   an expression
for the Hilbert function in terms
of the combinatorics of the Weyl group.
From this they recover the interpretation of the multiplicity,
due to Herzog and Trung \cite{ht}, as the cardinality of a certain set of
non-intersecting lattice paths.   One can interpret their main result as saying
that the natural set of generic determinantal minors that generate the ideal of
the tangent cone form a Gr\"obner basis in the ``anti-diagonal'' term order.
As to the case when $x$ is
some other point,  they conjecture expressions for
the Hilbert function and the multiplicity.

In the present paper,  we clarify the approach of Kreiman and Lakshmibai,
thereby handling, in their spirit,  all points.   We extend their
expressions
and combinatorial interpretation to all points.
Here again there is an interpretation of the main result in terms of Gr\"obner basis.

Even in the case when $x$ is the identity coset,  we have something to say
that is not already in \cite{kl}.     Namely we give explicitly a bijection between
the combinatorially defined sets $SM_w(m)$ and $S_w(m)$ of that paper.
The bulk of this paper,  namely \S4, is devoted to establishing 
this bijection.   While we are convinced that this bijection is natural and
moreover that it is, in some sense,  the only natural bijection,
we are unable to make this precise.   It would be good to have this
nailed down.

The expressions for the Hilbert function and
the multiplicity are given respectively in Theorem~\ref{tmain}  and
Corollary~\ref{cmain} of \S\ref{stheorem} below.
The combinatorial interpretation of
the multiplicity and the interpretation of the main result in terms of Gr\"obner basis
are given in \S\ref{sinterpretation}.
The proof of Theorem~\ref{tmain}
occupies \S\ref{ssmt} and \S\ref{sproof}.   In \S\ref{ssmt}
the proof of the Theorem~\ref{tmain} is reduced to the solution of a combinatorial problem,
and in \S\ref{sproof}
a solution of the combinatorial problem is presented.
It is only in \S\ref{sproof} that this paper differs significantly from
\cite{kl}.

Krattenthaler~\cite{k} has proved the multiplicity conjecture of \cite{kl} but there is
little in common between his approach and ours.    It is possible that the methods
of the recent paper of Knutson and Miller~\cite{km} are relevant to our situation but
we have not investigated further along this direction except to notice that the results
of \cite{kl}
can be deduced from those of \cite{km}.

The ideal of the tangent cone at the identity coset (as a sub-variety of
the tangent space of the Grassmannian) is a determinantal ideal,  the ideal 
``co-generated'' by a minor.    Abhyankar \cite{abh} and Herzog-Trung \cite{ht}
give formulas for the multiplicity and Hilbert function of the quotient
ring by such a determinantal ideal.   Lakshmibai-Weyman \cite{lw} obtain,
using standard monomial theory (for the Grassmannian, this theory boils down
to the results of Hodge
and Musili mentioned above),
a recursive formula for the
multiplicity at any point (not just at the identity coset).   From this
recursive formula Rosenthal-Zelevinsky \cite{rz}
obtain a closed-form determinantal
formula for the multiplicity at any point.  Krattenthaler \cite{kr1} gives
an interpretation of the Rosenthal-Zelevinsky formula as counting certain
non-intersecting lattice paths,  and then,    using the technique of
dual paths,  deduces a simpler determinantal formula.    
Krattenthaler's papers \cite{kr1} and \cite{k} provide
the bridge
between what we do here and the Rosenthal-Zelevinsky formula.

\section{The Theorem}\label{stheorem}
The purpose of this section is to state Theorem~\ref{tmain} and
Corollary~\ref{cmain} below which provide respectively combinatorial
expressions for  the Hilbert
function and the multiplicity at a point on a Schubert variety in
a Grassmannian.
We start by recalling some well known facts about Schubert varieties in
Grassmannians.     This also serves to fix notation.

Once for all fix positive integers $n$ and $d$ with $d\leq n$.
Let $\gdn$ denote the {\em Grassmannian},  which as a set consists of
$d$-dimensional linear subspaces of a fixed $n$-dimensional vector space.
Let us think of the $n$-dimensional vector space as the space of column matrices of
size $n\times 1$.   The general linear group $GL_n$ of $n\times n$ invertible
matrices acts on this by
left multiplication.   The induced action on $d$-dimensional subspaces
is transitive and so $\gdn$ is the quotient of $GL_n$ by the stabilizer
of any given point.   We identify $\gdn$
as the homogeneous
space $GL_n/P_d$,  where $P_d$ is the stabilizer of the point of $\gdn$
that is the span of the first $d$ standard basis vectors $e_1,\ldots,e_d$;
here $e_j$ is the $n\times 1$ column matrix whose entries are all zero
except the one on row $j$ which is $1$.

Let $B$ be the ``Borel subgroup'' of $GL_n$ consisting of
invertible upper triangular matrices,
and $T$ the ``torus'' of $GL_n$ consisting of invertible diagonal matrices.
The points of $\gdn$ that are fixed
by $T$ are precisely those $d$-planes that are spanned by $d$ of the standard
basis vectors. 
Let $\idn$ denote the set of subsets of cardinality $d$ of $\{1,\ldots,n\}$.
We write an element $v$ of $\idn$ as $(v_1,\ldots,v_d)$ where $1\leq v_1<\ldots
<v_d\leq n$.
Given an element $v=(v_1,\ldots,v_d)$ of $\idn$,
we denote by $e^v$ the $T$-fixed point of $\gdn$ that is the span of $e_{v_1},\ldots,
e_{v_d}$.
These points lie in different $B$-orbits and the union of the $B$-orbits
of these points is all of $\gdn$.

   A {\em Schubert variety} in $\gdn$
is by definition the closure of such a $B$-orbit (with the
reduced scheme structure).    It is a union of 
the $B$-orbit (of which it is the closure) with smaller dimensional
$B$-orbits.   Schubert varieties are thus indexed by $T$-fixed points
and so by $\idn$.    Given $w\in\idn$,  we denote by $X_w$ the closure of
the $B$-orbit of the $T$-fixed point $e^w$.

We are interested in the local rings of various points on a Schubert
variety $X_w$.    It is enough to focus attention on the $T$-fixed points
contained in $X_w$,  for $X_w$ is the union of $B$-orbits of such points.
The point $e^v$ belongs to $X_w$ if and only
if $v\leq w$,  where, 
  for elements $v= (v_1<\ldots<v_d)$ and $w=(w_1<\ldots<w_d)$
of $\idn$,   we say $v\leq w$ if $v_1\leq w_1$,
$\ldots$,
 $v_d\leq w_d$. 

For the rest of this section,  fix elements $v,w$ of $\idn$ with
$v\leq w$.  

Let $c$ be among the entries $v_1,\ldots,v_d$ of $v$, and
$r$ be among the complement $\{1,\ldots,n\}\setminus\{v_1,\ldots,v_d\}$.
Denote by $\roots^v$ the set of ordered pairs $(r,c)$ of such 
$r$ and $c$,  and by $\pos^v$ the subset of $\roots^v$ consisting
of those $(r,c)$ with $r>c$.  

We will be considering ``multisets'' of $\roots^v$ and $\pos^v$.
By a {\em multiset} we mean a subset in which repetitions are
allowed and kept account of.   Multisets can be thought of as
monomials in the variables consisting of elements of the set.
We also use the term  {\em monomial} to refer to a multiset.
The {\em degree} of a monomial is just the cardinality,
counting repetitions, of the multiset.

Given $\beta_1=(r_1,c_1)$, $\beta_2=(r_2,c_2)$ elements in $\pos^v$,
we say $\beta_1>\beta_2$ if $r_1>r_2$ and $c_1<c_2$.   A sequence
$\beta_1>\ldots>\beta_t$ of elements of $\pos^v$ is called a
{\em $v$-chain}.
Given a $v$-chain $\beta_1=(r_1,c_1)>\ldots>\beta_t=(r_t,c_t)$, we denote
by $s_{\beta_1}\cdots s_{\beta_t}v$ the element 
$(\{v_1,\ldots,v_d\}\setminus\{c_1,\ldots,c_t\})\cup\{r_1,\ldots,r_t\}$
of $\idn$.  In case the $v$-chain is empty,  this element is just $v$.
We say that $w$ {\em dominates} the $v$-chain
$\beta_1>\ldots>\beta_t$  if $w\geq 
s_{\beta_1}\cdots s_{\beta_t}v$.

Let $\mon$ be a monomial in $\roots^v$.  By a {\em $v$-chain in $\mon$}
we mean a sequence $\beta_1>\ldots>\beta_t$ of elements of
$\mon\cap\pos^v$.  We say that {\em $w$ dominates $\mon$} if $w$
dominates every $v$-chain in $\mon$.

Let $S^v_w$ denote the set of $w$-dominated monomials in $\roots^v$,
and $S^v_w(m)$ the set of such monomials of degree $m$.

We can now state our theorem:
\bthm \label{tmain} {\bf (Conjecture of Kreiman and Lakshmibai)}
Let $v,w$ be elements of $\idn$ with $v\leq w$.   Let $X_w$ be the
Schubert variety corresponding to $w$,  $e^v$ the $T$-fixed point
corresponding to $v$, and $R$ the co-ordinate ring of the tangent cone
to $X_w$ at the point $e^v$ (that is,  the associated graded ring
${\begin{rm}gr\end{rm}}_{\frak{M}}(O_{X_w, e^v})$ of the local ring
$O_{X_w,e^v}$ of $X_w$ at the point $e^v$ with respect to its maximal
ideal $\frak{M}$).     Then the vector space dimension of the 
$m^{\begin{rm}th\end{rm}}$ graded piece $R(m)$ of $R$ equals the 
cardinality of $\svwm$,  where $\svwm$ is as defined above.
\ethm
The proof of the theorem occupies sections \ref{ssmt} and \ref{sproof} .   
For now let us note the following easy consequence.
\bcor\label{cmain}
With notation as in Theorem~\ref{tmain} above,
the multiplicity of $R$ equals the number of square-free
$w$-dominated monomials in $\roots^v$ of maximum cardinality.
\ecor
\bproof
Consider the set $\left\{A_1,\ldots,A_k\right\}$ of all
members of $\svw$ that have no repetitions and that are maximal
with respect to inclusion among those having no repetitions. 
For a subset
$\ui:= \left\{i_1<\ldots<i_j\right\}$ of cardinality $j$ of
$\left\{1,\ldots,k\right\}$,  let $a_{\ui}$ denote the cardinality of
$A_{i_1}\cap\cdots\cap A_{i_j}$.   Then, since $\svwm$ is the set of
all monomials of degree $m$ with elements from any one of the $A_j$ as
variables,  we have
\[
	\begin{mbox}{\rm Cardinality of $\svwm$}\end{mbox}
=
\sum\limits_{j=1}^k (-1)^{j-1}\sum\limits_{\ui}
\binom{a_{\ui}-1+m}{a_{\ui}-1}
				\]
where the inner summation is over all subsets $\ui$ of cardinality $j$
of $\left\{1,\ldots,k\right\}$.   It follows immediately that the
multiplicity equals the number of $\ui$ with $a_{\ui}$ as large as
possible,  which  in other words is the number of $A_1,\ldots,A_k$ with
maximum cardinality.\eproof


\section{Reduction of the proof to combinatorics}\label{ssmt}
In this section we will see how some well known results of
Hodge~\cite[Chapter XIV]{hodge}
and Musili~\cite{smt} about homogeneous co-ordinate rings of Schubert
varieties in the Grassmannian allow us
to reduce the proof of Theorem~\ref{tmain} to the solution of 
a combinatorial problem.    This reduction is already there
in Lakshmibai-Weyman~\cite{lw} (see their Theorem~3.4) and in a
more general set-up,  but it would hurt readability if we merely
quoted their result and cut this section out.

Fix elements $v,w$ of $\idn$ with $v\leq w$, so that $e^v$ belongs to 
the Schubert variety $X_w$.   Let $\gdn\imbedding\Bbb{P}(\wedge^d V)$
be the Pl\"ucker embedding,  where we think of $\gdn$ as $d$-dimensional
subspaces of the $n$-dimensional vector space $V$.    For $\theta$
in $\idn$,  let $p_\theta$ denote the Pl\"ucker co-ordinate 
corresponding to $\theta$.    Let $\Bbb{A}^v$
be the affine patch of $\Bbb{P}(\wedge^d V)$ given by $p_v\neq0$,
 and set
$Y^v_w := X_w \intersect \Bbb{A}^v$.
The point $e^v$ is the origin of the affine space $\Bbb{A}^v$.

The functions $f_u := p_u/p_v$,  $u\in\idn$, provide a
set of co-ordinate functions on $\Bbb{A}^v$.  The co-ordinate ring
$k[Y^v_w]$ of $Y^v_w$ is a quotient of the polynomial ring 
$k[f_u\st u\in\idn]$,    $k$ being the underlying field. 
In order to describe a basis for $k[Y^v_w]$ as a $k$-vector space,   we now make some
definitions.

We consider multisets of $\idn$, or, in other words, monomials with elements
of $\idn$ as variables.
Such a monomial is called {\em standard}  if after
suitable rearrangement it has the form $\left\{\theta_1,\ldots,\theta_t\right\}$
with $\theta_1\geq\ldots\geq \theta_t$.   A standard monomial 
$\theta_1\geq\ldots\geq\theta_t$  is {\em $w$-dominated} if 
$w\geq\theta_1$;  it is {\em $v$-compatible}  if each of $\theta_1,
\ldots,\theta_t$ is comparable with $v$ and none of them equals $v$.
We denote the set of all $w$-dominated $v$-compatible standard monomials
by $\smvw$.   
\bprop\label{psmt} 
As $\theta_1\geq\ldots\geq\theta_t$ varies over $\smvw$, that is, over
all $w$-dominated $v$-compatible standard monomials,  the elements
$f_{\theta_1}\cdots f_{\theta_t}$ form a $k$-vector space basis of $k[Y_w^v]$.
\eprop
\bproof
We know from \cite[\S4]{smt} that the elements $p_{\tau_1}\cdots p_{\tau_s}$,
as $\tau_1\geq\ldots\geq\tau_s$  varies over all $w$-dominated standard
monomials,   form a basis for the homogeneous co-ordinate ring of $X_w$
as embedded in $\Bbb{P}(\wedge^d V)$.
(In \cite{smt}, $\gdn$ is identified as a right coset space rather than as a
left coset space as we have done here,  but there should be no problem in
translating the statements from there to our set up.)
Consider now any linear dependence relation among the $f_{\theta_1}
\cdots f_{\theta_t}$ with $\theta_1\geq \ldots\geq \theta_t$ in $\smvw$.
Multiplying this by a suitably high power of $p_v$ yields a linear
dependence relation among the $w$-dominated standard monomials
$p_{\tua_1}\cdots p_{\tau_s}$ (so called by an abuse of terminology),
and so the original dependence relation
must have only been the trivial one.


To prove that $f_{\theta_1}\cdots f_{\theta_t}$ generate $k[Y^v_w]$ as a vector
space (as $\theta_1\geq\ldots\geq\theta_t$ varies over $\smvw$), we use the
following fact:   if $\mu_1,\ldots,\mu_r$ be any monomial in $\idn$, and
$p_{\tau_1}\cdots p_{\tau_s}$ a standard monomial that occurs with non-zero
coefficient in the expression for $p_{\mu_1}\cdots p_{\mu_r}$ (as an element
of the homogeneous co-ordinate ring of $X_w$) as a linear combination of
$w$-dominated standard monomials,  then $r=s$ and $\tau_1\cup\cdots\cup\tau_s$
equals $\mu_1\cup\cdots\cup\mu_r$ as multisets of $\{1, 2,\ldots,n\}$.   This fact
follows from the nature of the quadratic relations (page 152 of \cite{smt}) which
provide the key to showing that the $w$-dominated standard monomials span
the homogeneous co-ordinate ring of $X_w$.
Let $\sigma_1,\ldots,\sigma_r$ be any monomial in
$\idn$.     Consider the  expression for
$p_{\sigma_1}\cdots
p_{\sigma_r}\cdot p_v^h$,  where $h$ is an integer such that $h>r(d-1)$, as a linear
combination of $w$-dominated standard monomials.
We claim that $p_v$ occurs in every standard monomial in this expression.
To prove the claim,  suppose that $p_{\tau_1}\cdots p_{\tau_{r+h}}$ is such a
standard monomial,  and that none of $\tau_1,\ldots,\tau_{r+h}$ equals $v$.  For
each $\tau_j$ there is at least one entry of $v$ that does not occur in $\tau_j$.
The number of occurrences of entries of $v$ in ${\tau_1}\cup\cdots\cup {\tau_{r+h}}$
is thus at most $(r+h)(d-1)$.   But these entries occur at least $hd$ times in
${\sigma_1}\cup\cdots\cup {\sigma_r}\cup v\cup\cdots\cup v$
(where $v$ is repeated $h$ times),  a contradiction to the fact above,
and the claim is proved.

Dividing by $p_v^{r+h}$ the expression for $p_{\sigma_1}\cdots p_{\sigma_r}\cdot
p^h_v$ as a linear combination of $w$-dominated standard monomials provides
an expression for $f_{\sigma_1}\cdots f_{\sigma_r}$ as a linear combination of
$f_{\theta_1}\cdots f_{\theta_t}$,   as $\theta_1\geq\ldots\geq\theta_t$ varies over
$\smvw$.\eproof

For $\theta\in\idn$,  define the {\em $v$-degree} of $\theta$
to be the cardinality of the set $\theta\setminus v$ (as a subset  of $\{1,2,\ldots,n\}$).
 Set $Z^v:=\gdn\cap\Bbb{A}^v$.  The co-ordinate ring $k[Z^v]$ of
$Z^v$ is a polynomial ring in the indeterminates $f_\theta$ with
$\theta$ varying over the elements of $\idn$ with $v$-degree $1$.
Such $\theta$ are in bijective correspondence with the elements
$\beta$ of $\roots^v$: if $\beta=(r,c)$,
then $\theta$ is obtained from
$v$ by replacing $c$ with $r$.   We set $X_\beta:=f_\theta$
and arrange them in a matrix as shown in Figure~1 
 for the case $d=4$, $n=10$, and $v=(2,5,7,9)$.
\[\begin{array}{c}
\left( \begin{array}{cccc}
	X_{12} & X_{15} & X_{17} &X_{19}\\
	1 & 0 & 0 & 0\\
	X_{32} & X_{35} & X_{37} &X_{39}\\
	X_{42} & X_{45} & X_{47} &X_{49}\\
	0 & 1 & 0 & 0\\
	X_{62} & X_{65} & X_{67} &X_{69}\\
	0 & 0 & 1 & 0\\
	X_{82} & X_{85} & X_{87} &X_{89}\\
	0 & 0 & 0 & 1\\
	X_{10,2} & X_{10,5} & X_{10,7} &X_{10,9}\\
\end{array} \right)\\
\mbox{}\\
\mbox{
Figure 1: The matrix of $X_\beta$, $\beta$ an element of $\roots^v$,}\\
\mbox{for the case $d=4$, $n=10$, and $v=(2,5,7,9).$}
\end{array}\]

The expression for a general $f_\theta$ in terms
of the indeterminates $X_\beta$ is obtained by taking the determinant
of the sub-matrix of such a matrix as in Figure~1 obtained by choosing
the rows given by the entries of $\theta$. 
Thus it is a graded polynomial of degree the $v$-degree of 
$\theta$. Since $X_w$ is the 
intersection of $\gdn$ with the planes $p_\tau$, $\tau\not\leq w$ (again see
\cite[Theorem~4.1]{smt}),
the co-ordinate ring $k[Y^v_w]$ is a quotient of the polynomial ring
$k[Z^v]$ by the homogeneous ideal $(f_\tau\st\tau\not\leq w)$.
We are interested in the tangent cone to $Y^v_w$ at the origin (the
point of $Z^v$ corresponding to $e^v$),  and since $k[Y^v_w]$
is graded,  this tangent cone is isomorphic to $k[Y^v_w]$ itself.

The proposition above tells us that the graded piece of degree $m$
of $k[Y^v_w]$ is generated as a $k$-vector space by elements of
$SM^v_w$ of degree $m$,  where the degree of a standard monomial
$f_{\theta_1}\cdots f_{\theta_t}$ is defined to be the sum of
the $v$-degrees of $\theta_1,\ldots,\theta_t$.
To prove Theorem~\ref{tmain} it therefore suffices
to show that the set $SM^v_w(m)$ of $v$-compatible $w$-dominated
standard monomials of degree $m$ is in bijection with $S^v_w(m)$.
\section{The Proof}\label{sproof}
In this section we show that $\smvw(m)$ and $\svw(m)$ are naturally
bijective.    As we saw in the previous section,   this serves to
complete the proof of Theorem~\ref{tmain}.

We fix once for all an element $v=(v_1,\ldots,v_d)$ of $\idn$.

We will describe two maps called $\pi$ and $\phi$.
The map $\pi$ associates
to any non-empty monomial $\mon$ of elements in $\pos^v$ a pair
$(w,\mon')$   consisting of an element $w$ of $\idn$ and a ``smaller''
monomial $\mon'$ in $\pos^v$.   This map enjoys the following good
properties:
\bprop\label{ppi}
\begin{enumerate}
\item
$w\gneq v$
\item
$\vdeg{w}+\degree{\mon'}=\degree{\mon}$
\item
$w$ dominates $\mon'$
\item
$w$ is the least element of $\idn$ that dominates $\mon$.
\end{enumerate}
\eprop
The map $\phi$, on the other hand,   associates a non-empty monomial in
$\pos^v$  to any pair $(w,\min)$ consisting of an element $w$ of $\idn$ with
$w\gneq v$  and a $w$-dominated monomial $\min$,  possibly empty, in $\pos^v$.
\bprop\label{ppiphi}
The maps $\pi$ and $\phi$ are inverses of each other.
\eprop

Before turning to the descriptions of $\pi$ and $\phi$ and the
proofs of the propositions above,
let us see how these mutually inverse maps help us establish the desired
bijection between $\smvw(m)$ and $\svw(m)$.
Let $S^v$ denote the set of monomials in $\roots^v$ and $T^v$ the set of monomials
in $\pos^v$.
Let $\smvtopv$ denote the set of $v$-compatible standard monomials that are anti-dominated
by $v$:   a standard monomial $\theta_1\geq \ldots\geq \theta_t$ is
{\em anti-dominated} by $v$ if $\theta_t\geq v$   (we can also write
$\theta_t>v$ since $\theta_t\neq v$ by $v$-compatibility).

Define the {\em domination map} from $T^v$ to $I(d,n)$ by
sending a monomial in $\pos^v$ to the least element that dominates it.
Define the {\em domination map} from $SM^{v,v}$
to $I(d,n)$ by sending $\theta_1\geq\ldots\geq\theta_t$ to
$\theta_1$. 
Both these maps take, by definition, the value $v$ on the empty monomial.

Using $\pi$, we now define a  map $\pitilde$  from $\tv$ to $\smvtopv$
that commutes with domination and preserves degree.    Proceed by
induction on the degree of an element $\mon$ of $\tv$.   
The image of the empty monomial under $\pitilde$ is taken to be the
empty monomial.
Let $\mon$ be non-empty,   and suppose that $\pi(\mon) =(w,\mon')$.    By (1) and
(2) of Proposition~\ref{ppi},   the degree of $\mon'$ is strictly less
than that of $\mon$,  and so by induction $\pitilde(\mon')$  is defined.
Suppose that $\pitilde(\mon')=w'\geq\ldots$.    By induction we also know
that the degree of $\mon'$ is the same as that of $w'\geq\ldots$ and that
$w'$ is the least element of $\idn$ that dominates $\mon'$.   By (3) of
Proposition~\ref{ppi},     we have $w\geq w'$,  and we set
$\pitilde (\mon):=w\geq\pitilde(\mon')$.   It follows from (4) of the same
proposition that $\mon$ has the same image under domination as
$\pitilde (\mon)$,  namely $w$.   It follows from (2) of the proposition
that the degrees of $\mon$ and $\pitilde(\mon)$ are the same, and the
induction is complete.

Using $\phi$, we  now define a map $\phitilde$ from $\smvtopv$ to $\tv$
that commutes with domination.    Proceed by induction on the length 
$t$ of a standard monomial  $\theta_1\geq\ldots\geq\theta_t$ in $\smvtopv$.
The image of the empty monomial under $\phitilde$ is taken to be the empty
monomial.   Letting $t\geq1$, we know by induction what $\phitilde(\theta_2\geq
\ldots\geq\theta_t)$ is,  and let us denote it by $\min$.   We also know
that $\theta_2$ is the least element that dominates $\min$, so that
$\theta_1$ dominates $\min$.    We set $\phitilde(\theta_1\geq\ldots\geq
\theta_t):=\phi(\theta_1,\min)$.    
By Proposition~\ref{ppiphi},  we know $\pi\phi(\theta_1,\min)=(\theta_1,\min)$,
so that the image under the domination map of $\pitilde\phi(\theta_1,\min)$ is
$\theta_1$.      Since $\pitilde$ commutes with domination,  it follows that $\theta_1$
is the least element to dominate $\phi(\theta_1,\min)$, and the induction is complete.

Easy arguments using induction,   which we omit,   show that $\pitilde$
and $\phitilde$ are mutually inverse maps.    Since $\pitilde$ preserves
degree,   it follows as a consequence that $\phitilde$ too preserves degree.

The maps $\pitilde$ and $\phitilde$   establish a bijection between
$\smvtopv$ and $\tv$.   In fact, since domination and degree are
respected,   we get a bijection $\smvtopvw(m)\cong\tvwm$, where
$\smvtopvwm$ is the set of $w$-dominated elements of $\smvtopv$
that have degree $m$,   and similarly $\tvwm$ is the set of $w$-dominated
elements of $\tv$ that have degree $m$.

Now let $U^v$ denote the set of monomials in $\roots^v\setminus\pos^v$
(that is,  in pairs $(r,c)$ of $\roots^v$ with $r<c$), and $SM^v_v$ 
the set of $v$-compatible and $v$-dominated standard monomials.  
As explained below,  the ``mirror image'' of the bijection
$SM^{v,v}\cong T^v$   gives a bijection $SM^v_v(m)\cong U^v(m)$,
where $SM^v_v(m)$ and $U^v(m)$ denote respectively the sets of elements
of degree $m$ of $SM^v_v$ and $U^v$.

Putting these bijections together, we get the desired
bijection:
\begin{eqnarray*}
\smvwm &=& 
	\bigcup\limits_{j=0}^m SM^{v,v}_w(j) \times
	SM^v_v(m-j)\\
&\cong&\bigcup\limits_{j=0}^m T^v_w(j) \times
				U^v(m-j)
				= S^v_w(m).
\end{eqnarray*}
Here the first equality is obtained by splitting a $v$-compatible
standard monomial $\theta_1\geq\ldots\geq\theta_t$ into two parts
$\theta_1\geq\ldots\geq\theta_p$ and $\theta_{p+1}\geq\ldots\geq\theta_t$,
where $p$ is the largest integer,  $1\leq p\leq t$,  with $\theta_p\geq v$
(we can also write $\theta_p>v$ since $\theta_p\neq v$ by
$v$-compatibility).    The last equality is obtained by writing a
monomial of degree $m$ in $\roots^v$ as a product of two monomials, one
in $\pos^v$ and the other in $\roots^v\setminus\pos^v$,
the sum of  their degrees being $m$.

Let us now briefly explain how to take the ``mirror image''.
For an integer $x$, $1\leq x\leq n$, define the dual $x^*$ by
$x^*:=n-x+1$.   For $v=(v_1,\ldots,v_d)$ in $\idn$,  define
the dual $v^*$ by $v^*:=(v_d^*,\ldots,v^*_1)$.   This dual map
on $\idn$ is an order reversing involution.  It induces a 
bijection $SM^v_v\cong SM^{v^*,v^*}$ by associating to
$\theta_1\geq\ldots\geq\theta_t$ the element $\theta_t^*
\geq\ldots\geq\theta^*_1$.   The sum of the $v$-degrees of 
$\theta_1,\ldots,\theta_t$ equals the sum of the $v^*$-degrees of
$\theta_t^*,\ldots,\theta_1^*$,  so that we get a bijection
$SM^v_v(m)\cong SM^{v^*,v^*}(m)$.

For an element $(r,c)$ in $\pos^{v^*}$,  define
the dual to be the element $(r^*,c^*)$
$\roots^v\setminus\pos^v$.
This induces a degree preserving bijection $T^{v^*}\cong U^v$.
Putting this together with the bijection of the previous paragraph and the
one given by $\pitilde$ and $\phitilde$ (for $v^*$ in place of $v$),
we have
\[
	SM^v_v(m)\cong SM^{v^*,v^*}(m)\cong T^{v^*}(m)\cong U^v(m)
								\]
The proof of Theorem~\ref{tmain} is thus reduced to the descriptions
of the maps $\pi$ and $\phi$ and the proofs of Propositions
\ref{ppi}  and \ref{ppiphi}.
\subsection{Preparation}\label{sspreparation}
The statements in the sequel and their proofs may appear
to be more involved than they really are.  We have found it convenient
to think of them in terms of pictures as illustrated in
Example~\ref{epi} below and in \S\ref{sinterpretation}.
In fact, they are just a translation into the language of
words of pictorial suggestions that are more or less obvious.

Recall from \S\ref{ssmt} the definition of $v$-degree of an element $\theta$ of
$\idn$:    it is the cardinality of the set $\theta\setminus v$ (or,
equivalently,  that of $v\setminus\theta$).   If $\beta_1>\ldots>\beta_t$ is a $v$-chain,
we call  $\beta_1$ the {\em head} of the $v$-chain, $\beta_t$ its {\em tail}, and $t$ its
{\em length}.
Let $\mon$ be a subset of $\pos^v$.    We say that an element $\beta$ of $\mon$
is {\em $t$-deep} in $\mon$ (where $t$ is a positive integer) if it is the tail of $v$-chain
in $\mon$ of length $t$.     The {\em depth} of $\beta$ in $\mon$ is $t$ if $\beta$ is
$t$-deep but not $(t+1)$-deep.   Two distinct elements $\alpha$ and $\beta$
of $\pos^v$ are said to be {\em comparable} if either $\alpha>\beta$ or $\beta>\alpha$.
Suppose, for example, that $\alpha=(r,c)$ and $\beta=(r,C)$.  In this case the two
are not comparable unless of course $c=C$.

The following proposition can be interpreted as specifying the map $\phi$
in the special case when $\min$ is empty.   It will be used in the description
of $\pi$.    Let us call {\em distinguished} the subsets (in the usual sense, not multisets) $\mon$
of $\pos^v$ satisfying the following conditions:
\begin{enumerate}
\item[(A)]
For $(r,c)\neq (r',c')$ in $\mon$,  we have $r\neq r'$ and $c\neq c'$
(in other words,  no two distinct elements of $\mon$ share a row or column
index).
\item[(B)]
If $\mon=\{(r_1,c_1),\ldots,(r_p,c_p)\}$ with $r_1<r_2<\ldots<r_p$,
then for $j$, $1\leq j\leq p-1$,  we have either $c_j>c_{j+1}$ or
$r_j<c_{j+1}$.
\end{enumerate}
Condition (B) can be restated as follows:
\begin{enumerate}
\item[(B*)]
For $(r,c), (R,C)$ in $\mon$ with $r<R$,  either $C<c$ or $r<C$.
\end{enumerate}
\bprop\label{psw}
There exists a bijection between elements $w$ of $\idn$ satisfying
$w\geq v$ on the one hand and subsets $\mon$ of $\pos^v$ satisfying
conditions (A) and (B) above on the other.  If $\mon_w$ be the subset
of $\pos^v$ corresponding to $w$ under this bijection, then
the $v$-degree of $w$ equals the cardinality of $\mon_w$.
\eprop
\bproof
Given $w\geq v$, consider the sets
$\{w_1,\ldots,w_d\}\setminus\{v_1,\ldots,v_d\}$
and
$\{v_1,\ldots,v_d\}\setminus\{w_1,\ldots,w_d\}$.
Both these have the same cardinality $p=\vdeg{w}$.   The first set provides us with
the row indices of elements of $\mon_w$,  the second with the column indices.
If we arrange the row indices in an increasing order,  say we have
$r_1<\cdots<r_p$,    then there is a unique way to arrange the
column indices such that condition (B*) above is satisfied: proceed by
induction,  and if $c_1,\ldots,c_j$ have been chosen,  choose $c_{j+1}$ to be
the maximum among the remaining column indices that are less than $r_{j+1}$.
This defines the map $w\mapsto\mon_w$.    It is clear that the cardinality
of $\mon_w$ equals $\vdeg{w}$.

For the map in the other direction,   given $\mon$,  to obtain $w$,  start with
$v=(v_1,\ldots,v_d)$,  delete those entries that occur as column indices in
$\mon$,  add those that occur as row indices in $\mon$,  and finally
arrange the entries in increasing order.  It is readily seen that the
two maps are inverses of each other.\eproof
\bremark\label{rsw}
\begin{enumerate}
\item
Subsets of distinguished monomials are themselves distinguished.
\item
If $\mon$ is a subset of $\pos^v$ satisfying condition (A),  then
we can still define a corresponding element $w$ of $\idn$ as in
the proof of the proposition.
If $\mon\subseteq\tilde{\mon}$ are subsets of $\pos^v$ satisfying condition
(A),   $w$ and $\tilde{w}$ the corresponding elements of $\idn$,  then
$w\leq \tilde{w}$.
\end{enumerate}
\eremark

\blemma\label{ldom}
Let $\beta_1=(r_1,c_1)>\ldots>\beta_t=(r_t,c_t)$ be a $v$-chain,
$w$ and element of $\idn$ with $w\geq v$, and $\monw$ the
distinguished subset of $\pos^v$ associated to $w$.  Then
$w$ dominates $\beta_1>\ldots>\beta_t$ if and only if there
exists a $v$-chain $\alpha_1=(R_1,C_1)>\ldots>\alpha_t=(R_t,C_t)$
in $\monw$ such that $C_j\leq c_j$ and $r_j\leq R_j$
for $1\leq j\leq t$.
\elemma
\bprooftwo
Just for the purpose of this proof,  we introduce the
following notation:   for a subset $A$ of positive integers
and $l$ a positive integer,  denote by $A^l$ the subset
$\{a\in A\st a\leq l\}$.   Note that for elements $u,y$ of
$\idn$,  we have $u\leq y$ if and only if $\card{u^l}
\geq \card{y^l}$ for every $l$,  where $\card{\cdot}$
denotes cardinality.

Suppose first that there exists a $v$-chain $\alpha_1>\ldots>\alpha_t$
of elements of $\monw$ as in the statement.   Set
$y:=s_{\alpha_1}\cdots s_{\alpha_t}v$ and
$u:=s_{\beta_1}\cdots s_{\beta_t}v$.   Clearly $w\geq y$
(see Remark~\ref{rsw}~(2)).   To show that $w$ dominates
$\beta_1>\ldots>\beta_t$ it therefore suffices to show $y\geq u$.

Let $l$ be any integer.    We have
\begin{eqnarray*}
u^l&=&\left(v^l\setminus\left\{c_1,\ldots,c_t\right\}^l\right)
		\cup\left\{r_1,\ldots,r_t\right\}^l
		\quad\quad \mbox{\begin{rm} and \end{rm}}
		\quad\quad \\
y^l&=&\left(v^l\setminus\left\{C_1,\ldots,C_t\right\}^l\right)
		\cup\left\{R_1,\ldots,R_t\right\}^l
						\end{eqnarray*}
From the hypothesis it immediately follows that
\[
	\left|\left\{c_1\ldots,c_t\right\}^l\right|\leq
	\left|\left\{C_1\ldots,C_t\right\}^l\right|
		\quad\quad \mbox{\begin{rm} and \end{rm}}
		\quad\quad
	\left|\left\{r_1\ldots,r_t\right\}^l\right|\geq
	\left|\left\{R_1\ldots,R_t\right\}^l\right|
						\]
Thus $\card{u^l}\geq\card{y^l}$  and we are done
(implicitly used here are the facts that $\{c_1,\ldots,c_t\}$,
$\{C_1,\ldots,C_t\}$ are subsets of $v$ and $\{r_1,\ldots,r_t\}$,
$\{R_1\ldots,R_t\}$ are disjoint from $v$).

To prove the converse,  we proceed by induction on $t$.
Let us isolate the case $t=1$ as a separate statement:
\blemma\label{lsubldom}
Let $\beta=(r,c)$ be an element of $\pos^v$ and $w$ an element
of $\idn$ such that $w\geq s_\beta v$.   Then there exists
an element $(R,C)$ in the distinguished monomial $\monw$
associated to $w$ such that $C\leq c$ and $r\leq R$.
\elemma
\bproof
Consider the set of elements $(A,B)$ of $\monw$ such that
$B\leq c<A$.   We claim that this set is non-empty.  To see
this,   consider $w^c$ and $u^c$,  where $u:=s_\beta v$.
Clearly $u^c=v^c\setminus\{c\}$.   If the set were empty,
we would have $\card{w^c}=\card{v^c}$, and so $\card{w^c}
\gneq\card{u^c}$,  a contradiction to the hypothesis that
$w\geq u$.

Since $\monw$ is distinguished,  the set of such elements
$(A,B)$ form a $v$-chain.     Let $(R,C)$ be the head of this
$v$-chain.     We claim that $r\leq R$.    To see this,
consider $w^R$ and $u^R$.     Since $\monw$ is distinguished,
we have $\card{w^R}=\card{v^R}$.   If $R<r$,  then
$u^R=v^R\setminus\{c\}$,   and so $\card{w^R}\gneq\card{v^R}$,
a contradiction to the hypothesis that $w\geq u$.\eproof

Now suppose that $w$ dominates $\beta_1>\ldots>\beta_t$.
Then $w\geq s_{\beta_1}\cdots s_{\beta_t}v\geq s_{\beta_1}v$
(see Remark~\ref{rsw}).    By the lemma above,   there exists
$(R,C)$ in $\monw$ with $C\leq c_1$ and $r_1\leq R$.  Since
$\monw$ is distinguished,  the set of such $(R,C)$ forms a
$v$-chain.    Let $\alpha_1=(R_1,C_1)$ be the head of this
$v$-chain.

Let $\tilde{w}$ be the element of $\idn$ associated to the
distinguished subset $\monw\setminus\{(R_1,C_1)\}$ of $\pos^v$.
We will show below that $\tilde{w}$ dominates $\beta_2>\ldots>\beta_t$.
Assuming this for the moment,  we obtain by induction on $t$
(the case $t=1$ being the lemma above) a $v$-chain $\alpha_2=
(R_2,C_2)>\ldots>\alpha_t=(R_t,C_t)$ with $C_j\leq c_j$ and $r_j\leq
R_j$.  Since $\monw$ is distinguished,   it follows easily from
our choice of $\alpha_1$ that $\alpha_1>\alpha_2$, and we are done.

To prove that $\tilde{w}\geq s_{\beta_2}\cdots s_{\beta_t}v$, we
set $y:=s_{\beta_2}\cdots s_{\beta_t}v$ and verify that, for
every integer $l$, $\card{\tilde{w}^l}\leq\card{y^l}$:
\begin{itemize}
\item
	For $l<c_2$, we have $y^l=v^l$.  Since $\tilde{w}\geq v$,
	it follows that $\card{\tilde{w}^l}\leq\card{y^l}$.
\item
	For $c_2\leq l<r_1$,  we have $y^l=u^l\cup\{c_1\}$,
	so that $\card{y^l}=\card{u^l}+1$,  and $\tilde{w}^l
	=w^l\cup\{C_1\}$.  Since $w\geq u$, it follows that
	$\card{\tilde{w}^l}=\card{w^l}+1\leq\card{u^l}+1=\card{y^l}$.
\item
	For $l\geq r_2$,  we have $\card{y^l}=\card{v^l}$. Since
		$\tilde{w}\geq v$,  it follows that $\card{
			\tilde{w}^l}\leq\card{y^l}$.~$\Box$
			\end{itemize}\end{prooftwo}

\bcor\label{cldom}
Let $w$ be an element of $\idn$ and $\monw$ the corresponding
distinguished subset of $\pos^v$.  For a positive integer $j$,
let $\mon_w^j$ denote the subset of $\monw$ of those elements
that are $j$-deep,  and $w^j$ the corresponding element of
$\idn$.  Let $\beta_1=(r_1,c_1)>\ldots>\beta_t=(r_t,c_t)$ be
a $v$-chain.
\begin{enumerate}
\item
	If $w^k$ dominates $\beta_1>\ldots>\beta_t$,  then $w^{k+1}$
	dominates $\beta_2>\ldots>\beta_t$, $w^{k+2}$ dominates
	$\beta_3>\ldots>\beta_t$, and so on.
\item
                If, for integers $l>k$,  there exists $(R,C)$ in $\monw^l$ such
                that $C\leq c_1$ and $r_1\leq R$,  and $w^{k+1}$ does not dominate
                $\beta_1>\dots>\beta_t$,  then $w^{k+2}$ does not dominate
               $\beta_2>\dots>\beta_t$,   $w^{k+3}$ does not dominate
                $\beta_3>\dots>\beta_t$,   and so on until, finally,
                $w^{l+1}$ does not dominate
                $\beta_{l-k+1}>\dots>\beta_t$.
\end{enumerate}
\ecor
\bproof
(1)~~Choose a $v$-chain $\alpha_1>\ldots>\alpha_t$ in $\mon_w^k$
as in the lemma.  Note that $\{\alpha_2,\ldots,\alpha_t\}
\subseteq\mon_w^{k+1}$,  $\{\alpha_3,\ldots,\alpha_t\}\subseteq
\mon_w^{k+2}$, and so on.   Now use the ``if'' part of the
lemma.

(2)~~The proof is based on the following claim:  given an element
$\beta$ of $\pos^v$ and an integer $j$,  there exists at most one
element of $\alpha$ of $\mon_{w,j}$ such that $\alpha>\beta$.  To prove
the claim,   note that, on the one hand,  since $\monw$ is distinguished,
the set of its elements $\alpha$ such that $\alpha>\beta$ form a $v$-chain,
but, on the other hand,   no two elements of $\mon_{w,j}$ are comparable.

It is enough to do the case when $l=k+1$.  So let $l=k+1$.
Suppose that $w^{k+2}$ dominates $\beta_2>\ldots>\beta_t$.
By the lemma above,  there exists a $v$-chain $(R_2,C_2)>\ldots
>(R_t,C_t)$ in $\mon_w^{k+2}$ such that $C_j\leq c_j$ and $r_j\leq R_j$
for $2\leq j\leq t$.     Choose $(R_1,C_1)$ in $\mon_{w,k+1}$ such that
$(R_1,C_1)>(R_2,C_2)$.  If $(R,C)$ belongs to $\mon_w^{k+2}$,  then
replacing it with an element $(R',C')$ of $\mon_{w,k+1}$ such that
$(R',C')>(R,C)$,  we may assume that $(R,C)$ belongs to $\mon_{w,k+1}$.

Now we have $(R_1,C_1)>\beta_2$ and $(R,C)>\beta_2$.  It follows from
the claim above that $(R_1,C_1)=(R,C)$.  Thus $C_1\leq c_1$ and
$r_1\leq R_1$, and, by the ``if'' part of the lemma above,  it follows
that $w^{k+1}$ dominates $\beta_1>\ldots>\beta_t$.\eproof

\subsection{The description of $\pi$}\label{sspi}
In what follows,
when we speak of a ``subset'' $\mon_1$ of a
multiset $\mon$,  it should be understood that we have in mind a subset of the
set underlying $\mon$ and that $\mon_1$ is the collection
with multiplicities of those elements of $\mon$ that belong to the
underlying subset.    A similar comment
applies to ``partition'' of a multiset.

We now specify the map $\pi$.
Let $\mon$ be a non-empty monomial in the elements of $\pos^v$.
We partition $\mon$
in two stages.  First we partition $\mon$ into subsets
$\mon_1,\ldots,\mon_k$,   where $k$ is the largest length of
a $v$-chain in $\mon$:   $\beta\in\mon$ belongs to $\mon_j$
if it is $j$-deep but not $(j+1)$-deep.
Clearly $\mon_1,\ldots,\mon_k$ form a partition of $\mon$,
and no two distinct elements belonging to the same
$\mon_j$ are comparable.

Now we partition each $\mon_j$ into subsets called
{\em blocks} as follows.   We arrange the elements of $\mon_j$
in non-decreasing order of their row numbers (all arrangements are from left
to right; and elements occur with their respective multiplicities).
Among those with
the  same row number,  the arrangement is by non-decreasing
order of column numbers.   By the incomparability of any two
distinct elements of $\mon_j$,    the column numbers are
also now arranged in non-decreasing order.   If $(r,c),(R,C)$
are consecutive members in this arrangement,   we demand that
the two belong to the same or different blocks accordingly as
$r>C$ or $r<C$ ($r\neq C$, for $C$ belongs and $r$ does not to
$\{v_1, \ldots, v_d\}$).    More formally,  two consecutive
members $(r,c),(R,C)$ in the arrangement are said to be {\em related}
if $r>C$.    The blocks are the equivalence classes
of the smallest equivalence relation containing the above relations.

Let us now concentrate on a single block $\part$ of some $\mon_j$.
Let \begin{eqnarray}\label{eblock}
	(r_1,c_1), \ldots, (r_p,c_p)	\end{eqnarray}
be the elements of $\part$ written in non-decreasing order of
both row and column numbers (in all such arrangements
the elements occur with their
respective multiplicities).    We set $w(\part):= (r_p,c_1)$,
$\part'$ to be the monomial
\begin{eqnarray} \label{eblockprime}
		\left\{(r_1,c_2), (r_2,c_3), \ldots,
		(r_{p-2},c_{p-1}), (r_{p-1}, c_p)\right\} \\
\label{emonprime}
\mon_j':=\bigcup\limits_{\block} \block'   \quad\quad
\mbox{\begin{rm} and \end{rm}}\quad\quad
\mon':=
 \bigcup\limits_{j=1}^k \mon_j'
			\end{eqnarray}
where the index $\part$ runs over all
blocks of $\mon_j$. 
%
It follows from Corollary~\ref{cltwoone} below that
\[      \left\{w(\block)\st\mbox{\begin{rm} $\block$ is a block of $\mon$ \end{rm}}
                                                                                                                \right\}\]
satisfies conditions (A) and (B) of subsection~\ref{sspreparation}.  Let $w$
be the corresponding element of $\idn$ (see Proposition~\ref{psw}),  and
set
\[	\pi(\mon) := (w,\mon').		\]

\renewcommand{\thefigure}{\arabic{figure}}
Let us illustrate the definition of the map $\pi$ by means of
an example.
\bexample\label{epi}
Let $d=13$, $n\geq 25$,  and
$$v=(1,2,4,5,7,8,9,14,15,16,17,18,19).$$
A monomial $\mon$ and its decomposition into blocks is represented
in Figure 2.

\begin{figure}\label{fexmon}
\begin{center}
\mbox{\epsfig{file=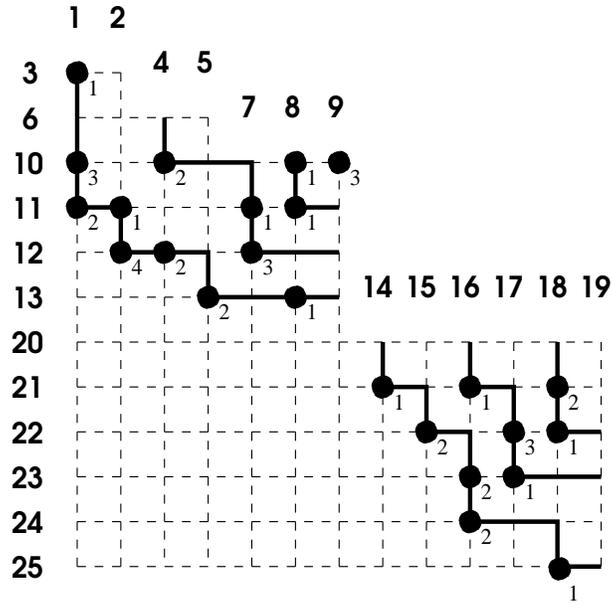,height=8cm,width=8cm}}
\caption{A monomial $\mon$ and its decomposition into blocks}
\end{center}
\end{figure}

The lattice points of the grid correspond to elements of $\pos^v$.
The solid dots indicate the elements that
occur with non-zero multiplicity.   The integers by the dots are
the multiplicities of the corresponding elements.   The solid lines
indicate the decomposition of $\mon$ into blocks.   There are two
blocks each in $\mon_1$, $\mon_2$, and $\mon_3$.  For example, the
blocks in $\mon_3$ are $\{(10,8),(11,8)\}$ and $\{(21,18)^2, (22,18)\}$
(the exponent $2$ indicates the multiplicity);   $\mon_4$ consists
of just one block, namely, $\{(10,9)^3\}$.

Let $\pi(\mon)=(w,\mon')$.  Since $\mon$ has  $7$ blocks,
$\mon_w$ has $7$ elements:
\[
	\mon_w=\left\{(10,9), (11,8), (12,4), (13,1),
 (22,18), (23,16), (25,14)\right\}		\]
The monomial $\mon'$ is represented in
Figure~3.\eexample

\begin{figure}\label{fexmonprime}
\begin{center}
\mbox{\epsfig{file=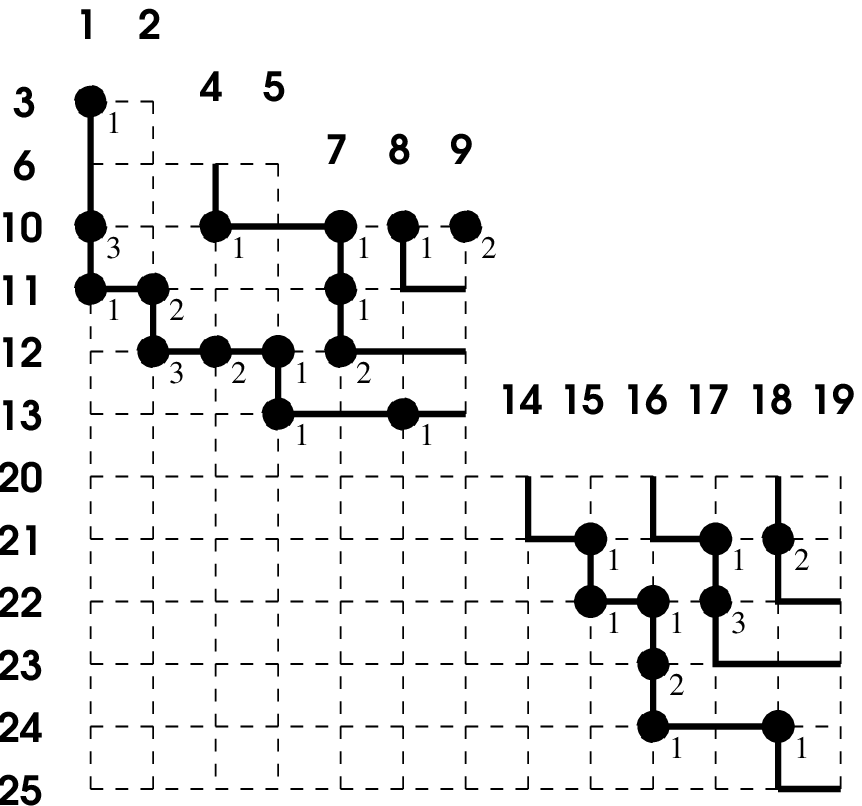,height=8cm,width=8cm}}
\caption{The monomial $\mon^\prime$ where $\mon$ is the monomial of Figure~2}
\end{center}
\end{figure}

We now note some properties of the above construction that will be
used in the proofs of Propositions \ref{ppi} and \ref{ppiphi}.
\blem\label{lone}
If $\block_1,\ldots, \block_l$ are the blocks in order
 from left to right of  some $\mon_j$,   and
$w(\block_1)=(R_1,C_1)$,
$w(\block_2)=(R_2,C_2)$,
$\ldots$,
$w(\block_l)=(R_l,C_l)$,
then \[ C_1<R_1<C_2<R_2<\ldots<R_{l-1}<C_l<R_l	\]
\elem
\bproof
From the definition of $w(\block)$,  it is clear that
$C_1<R_1$,  $C_2<R_2$, etc.   Note that $C_2$ is the least
column index in $\block_2$ and $R_1$ is the greatest row
index in $\block_1$.   From the way we divided $\mon_j$
into blocks,  it is clear that $R_1<C_2$.  Similarly
$R_2<C_3$, etc.\eproof
\blemma\label{lssprime}
No two
distinct elements of $\mon_j\cup\mon_j'$ are comparable.
More precisely,   if $(a,b)$ and $(c,d)$ are elements of $\pos^v$
with $(a,b)>(c,d)$,   then they cannot both belong to
$\mon_j\cup\mon_j'$.
\elemma
\bproof
We can put the elements in (\ref{eblock}) and (\ref{eblockprime})
together in non-decreasing order of both row and column numbers:
\begin{eqnarray*}
(r_1,c_1), (r_1,c_2), (r_2,c_2), (r_2,c_3), (r_3,c_3),\ldots,\\
(r_{p-2},c_{p-2}), (r_{p-2},c_{p-1}), (r_{p-1},c_{p-1}),
	(r_{p-1},c_p),(r_p,c_p) 		\end{eqnarray*}
This means that distinct  elements belonging to $\block\cup\block'$,
where $\block$ is a block of $\mon_j$,  are not comparable.
Let $\block$ and $\blockc$ are distinct blocks of $\mon_j$
with $w(\block)=(R,C)$ and $w(\blockc)=(S,T)$.  Using
Lemma~\ref{lone}, we may assume without loss of generality that
$C<R<S<T$.    If $(a,b)$ belongs to $\block\cup\block'$ and
$(c,d)$ to $\blockc\cup\blockc'$,  then
\[	C\leq b< a\leq R<S\leq d<c\leq T	\]
so that $(a,b)$ and $(c,d)$ are not comparable.
Note that this shows also that $\block\cup\block'$
and $\blockc\cup\blockc'$ are disjoint.\eproof
\bcor\label{clssprime}
For $(r,c)\in\mon'$,  there exists a unique block $\block$ of
$\mon$ with $(r,c)\in\block'$.
\ecor
\bproof
The existence is immediate from the definition of $\mon'$.
For the uniqueness,   let $\block$ and $\blockc$ be two distinct
blocks of $\mon$ with $(r,c)$ belonging to both $\block'$ and
$\blockc'$.    For the reason noted at the end of the proof of
the previous lemma,   $\block$ and $\blockc$ cannot belong to
the same $\mon_j$.     So suppose, without loss of generality,
that $\block\subseteq\mon_i$, $\blockc\subseteq\mon_j$, and $i<j$.
Now there exists an element of the form $(r,a)$ in $\blockc$ with
$a\leq c$ (this follows from the definition of $\blockc'$),  and
an $(R,C)$ in $\mon_i$ with $(R,C)>(r,a)$ (such an
element can be found at the appropriate place on a $v$-chain of
length $j$ with $(r,a)$ as tail).    But then $(R,C)>(r,c)$
and both $(R,C)$ and $(r,c)$ belong to $\mon_i\cup\mon_i'$,
a contradiction to the previous lemma.\eproof
\blem\label{ltwo}
Given positive integers $i<j$ and a block $\block$ of $\mon_j$,
there exists a unique block $\blockc$ of $\mon_i$ such that
$w(\blockc)>w(\block)$.
\elem
\bproof
Let $\blc$ and $\blcone$ be distinct blocks of $\mon_i$
with $w(\blc)=(A,B)$, $w(\blcone)=(A_1,B_1)$, and
$(A,B)>(r,c)$.     By Lemma~\ref{lone},
\[\mbox{\begin{rm}
either
 $B<c<r<A<B_1<A_1$ or
 $B_1<A_1<B<c<r<A$.\end{rm}}\]
In either case, we have $(A_1,B_1)\not>(r,c)$.  This proves uniqueness.

For the existence,  write $w(\block)=(r,c)$, and
let the elements of $\block$ in non-decreasing
order of row and column entries be
\[  (r_1,c),\ldots,(r,c_p)	\]
Let $(R,C)$ and $(R',C')$  be elements of $\mon_i$ such that
$(R,C) > (r_1,c)$ and $(R',C')>(r,c_p)$.
If $(R',C')$ either equals $(R,C)$ or appears to the left of it
in the arrangement of the elements of $\mon_i$ in non-decreasing
order of row and column entries,  then let $\blockc$ be the
block of $\mon_i$ containing $(R',C')$.    Writing $w(\blockc)
=(A,B)$,  we have $r<R'\leq A$ and $B\leq C'\leq C<c$,
and we are done.

So let us assume
that $(R',C')$ appears to the right of $(R,C)$  in the arrangement
of elements of $\mon_i$:
\[	\ldots,(R,C),\ldots,(s,t),\ldots,(R',C'),\ldots	\]
where $(s,t)$ denotes an arbitrary element between $(R,C)$ and $(R',C')$ that is not equal to
$(R',C')$.    We claim that $s>b$ for an element
$(a,b)$ that occurs in the above arrangement
strictly to the right of $(s,t)$.
Once the claim is proved,  it will follow that
all the elements $(R,C),\ldots,(R',C')$
belong to the same block of $\mon_i$.    Let $\blc$ be this
block and set $w(\blc)=(A,B)$.    We then have
$B\leq C<c$ and $r<R'\leq A$.

%
%
To prove the claim,  we may assume that $s<c_p$, for otherwise $C'<c_p\leq s$
and we can take $(a,b)$ to be $(R',C')$.   Let $y$ be the least
among the column indices of elements of $\block$ subject to
$s<y$,  and let $(x,y)$ be the left most element with column
index $y$ in the arrangement of $\block$.   It is not possible
that $(r_1,c)=(x,y)$,  for that would imply $y<x=r_1<R\leq s$,
but $s<y$ by our choice of $y$.

Now let $(x',y')$ be the element
in $\block$ immediately to the left of $(x,y)$.
By our choice of $(x,y)$ we have $y'\leq s$.  Since $(x',y')$ and
$(x,y)$ belong to the same block $\block$,  we have $y<x'$.

Let $(a,b)$ be an element in $\mon_i$ such that $(a,b)>(x',y')$.
We then have $s<y<x'<a$, so that $(a,b)$ appears to the right of
$(s,t)$.  Further,  we have $b<y'\leq s$, and the claim is proved.\eproof
\bcor\label{cltwoone}
If $\block$ and $\block_1$ are blocks of $\mon$ with $w(\block)=(r,c)$
and $w(\block_1)=(r_1,c_1)$,  then exactly one of the following holds:
\begin{eqnarray*}
        c<r<c_1<r_1,&\quad\quad\quad\quad& c_1<r_1<c<r,\\
  c<c_1<r_1<r,&\quad\quad\mbox{\begin{rm}or\end{rm}}
  \quad\quad& c_1<c<r<r_1.
  \end{eqnarray*}
  \ecor
  \bproof
  Let $j$ and $j_1$ be the integers such that $\block\subseteq\mon_j$
  and $\block_1\subseteq\mon_{j_1}$.   If $j=j_1$,  then it follows from
  Lemma~\ref{lone} that one of the first two possibilities holds.

  Now suppose that $j_1<j$.   Then,  by Lemma~\ref{ltwo},   there exists
  a block $\blockc$ of $\mon_{j_1}$ such that $w(\blockc)>w(\block)$.
  Writing $w(\blockc)=(A,B)$,  we have $B<c<r<A$.   If $\block_1=\blockc$,
  then the fourth possibility holds.   If $\block_1\neq \blockc$,  then, again
  by Lemma~\ref{lone},  either $B<A<c_1<r_1$ or $c_1<r_1<B<A$,  so
  that respectively either the first or the second possibility holds.

  If $j_1<j$,   then we see similarly that one of the first three possibilities
  holds.\eproof
\bcor\label{cltwo}
\begin{enumerate}
\item
Let $(r,c)>(r_1,c_1)$ be elements of $\mon$, and $\block$,
$\block_1$ be blocks of $\mon$ such that $(r,c)\in\block$,
and $(r_1,c_1)\in\block_1$.  Then $w(\block)>w(\block_1)$.
\item
Let $(r,c)>(r_1,c_1)$ be elements of $\mon'$, and $\block$,
$\block_1$ be blocks of $\mon$ such that $(r,c)\in\block'$,
and $(r_1,c_1)\in\block_1'$.  Then $w(\block)>w(\block_1)$.
\end{enumerate}
\ecor
\bproof  We prove both statements simultaneously.
Let $j$ and $j_1$ be  such that $\block\subseteq\mon_j$
and $\block_1\subseteq\mon_{j_1}$.  We claim that $j<j_1$.
This is clear in case (1).    In case (2),
$j\neq j_1$ by
Lemma~\ref{lssprime}.   If $j>j_1$,
then choosing an element of the form $(r,a)$ with
$a\leq c$ in $\block$  and an element $(R,C)$ in $\mon_{j_1}$ with $(R,C)>(r,a)$, we have
$(R,C)>(r_1,c_1)$,  but this contradicts Lemma~\ref{lssprime}.
The claim is thus proved.

Write $w(\block)=(R,C)$ and 
$w(\block_1)=(R_1,C_1)$.   By Lemma~\ref{ltwo},  there exists
a unique block $\blockc$ of $\mon_j$ with
$w(\blockc)>w(\block_1)=(R_1,C_1)$.  We will show that
$\block=\blockc$.
Write $w(\blockc)=(S,T)$ and suppose that $\block\neq\blockc$.
Then, by Lemma~\ref{lone},  either  $S<T<C<R$ or $C<R<S<T$.
In the first case we have
\[
	S<C_1\leq c_1<r_1\leq  R_1<T<C\leq c<r\leq R	\]
	which is inconsistent with our hypothesis that
		$c<c_1<r_1<r$.    A similar contradiction occurs
			in the second case.\eproof  
\bcor\label{cltwothree}
Let $\beta$ be an element of $\mon$,  $\block$ the block of 
$\mon$ containing $\beta$,  and $\block_1$ a block of $\mon$
with $w(\block_1)>w(\block)$.  Then there exists $\beta_1$
in $\block_1$ with $\beta_1>\beta$.
\ecor
\bproof
Let $j$ and $j_1$ be the integers such that $\block\subseteq
\mon_j$ and $\block_1\subseteq\mon_{j_1}$.  By Lemma~\ref{lone},
$j\neq j_1$.  Suppose $j<j_1$.  By Lemma~\ref{ltwo},  there exists
a block $\blockc$ in $\mon_j$ with $w(\blockc)>w(\block_1)$,
but then $w(\blockc)>w(\block)$,  which contradicts Lemma~\ref{lone}.
Thus $j>j_1$.

Now there exists $\beta_1$ in $\mon_{j_1}$ such that $\beta_1>\beta$.
Let $\blockc$ be the block of $\mon_{j_1}$ containing $\beta_1$.
By Corollary~\ref{cltwo}~(1),  $w(\blockc)>w(\block)$,
and by the uniqueness part of Lemma~\ref{ltwo}, 
$\blockc=\block_1$.\eproof
\bcor\label{cltwofour}
Let $\block$ be block of $\mon$ and $\beta$ an element of
$\block$.  Then the depth of $\beta$ in $\mon$ equals the
depth of $w(\block)$ in $\left\{w(\block)\st\mbox{
\begin{rm}$\block$ is a block of $\mon$\end{rm}}\right\}$.
\ecor
\bproof
That the depth of $w(\block)$ is not less than the depth of
$\beta$ follows from Corollary~\ref{cltwo}~(1).    That the
depth of $\beta$ is not less than that of $w(\block)$
follows from Corollary~\ref{cltwothree}.\eproof

\subsection{Proof of Proposition~\ref{ppi}} \label{sspiproof}
It is clear from the definition of $\mon'$ that $\degree{\mon'}$ is less
than $\degree{\mon}$ by exactly the number of blocks in $\mon$.
By Proposition~\ref{psw},  this number equals $\vdeg{w}$,   and so
statement (2) is proved.   Moreover,   this number is positive since
$\mon$ is non-empty,  and so statement (1) holds.


	To prove (3), let a $v$-chain $\beta_1>\ldots>\beta_t$
	in $\mon'$ be given.   Let $\block_1,\ldots,\block_t$
	be blocks of $\mon$ such that $\beta_i$ belongs to $\block_i'$.
	(In fact,  these blocks are uniquely determined as
	shown in Lemma~\ref{clssprime},  but here we need only
	that they exist, and this is immediate from the definition
	of $\mon'$.)    It follows from Corollary~\ref{cltwo}~(2) that
	$w(\block_1)>\ldots>w(\block_t)$,  and then from 
	Lemma~\ref{ldom} that (3) holds.

%

	We now turn to the proof of (4).   That $w$ dominates
	$\mon$ follows in much the same way as statement (3):
                it uses Corollary~\ref{cltwo}~(1) in place of Corollary~\ref{cltwo}~(2).
               To prove that it is
	the minimal one that does so,  it suffices to show
	the following:  given an integer $j$, $1\leq j\leq d$,  there
	exists a $v$-chain $\alpha_1>\ldots>\alpha_t$ in $\mon$
	such that $y:=s_{\alpha_1}\cdots s_{\alpha_t}v$ satisfies
	$y_j=w_j$,  where $y_j$ and $w_j$ denote as usual the
	$j^{\mbox{\begin{rm}th\end{rm}}}$ entries of $y$ and
	$w$ as elements of $\idn$.

	Fix $j$, $1\leq j\leq d$.  We will first get hold of
	a $v$-chain $\beta_1>\ldots>\beta_t$ in
        \(\mon_w:=\left\{w(\block)\st
        \mbox{\begin{rm} $\block$ is a block
	of $\mon$ \end{rm}}\right\}\)
	such that $u=s_{\beta_1}\cdots s_{\beta_t}v$ satisfies
	$u_j=w_j$.   Set $D:=w_j$.  Consider the subset of
	$\mon_w$ consisting of those elements $(r,c)$ for which
	$c<D\leq r$.
                Since $\mon_w$ is distinguished
                (see Corollary~\ref{cltwoone}),
                it follows that
	these elements form a $v$-chain,  say $\beta_1=(r_1,c_1)>
	\ldots>\beta_t=(r_t,c_t)$.    We have
                \begin{eqnarray}
                        c_1<\ldots<c_t<D\leq r_t<\ldots<r_1    \label{epiproofone} \end{eqnarray}
                        and $r_t=D$ if and only if $D$ is not an entry of $v$.

	Set $u:=s_{\beta_1}\cdots
	s_{\beta_t}v$.   We claim that $u_j=D=w_j$.  First note
	that $D$ occurs as an entry of  $u$:  if $D$ belongs to
	$\{v_1,\ldots,v_d\}$, then $D$ is not a column index
	in $\mon_w$  much less of any $\beta_1,\ldots,\beta_t$;
	if $D$ is a row index in $\mon_w$,  then $r_t=D$.
	Every element $\delta=(r,c)$
	of $\mon_w$ other than $\beta_1,\ldots,\beta_t$
	satisfies either $c<r<D$ or $D<c<r$,  so that $s_\delta u$
	will continue to have $D$ in the same place as $u$.  Since
	we can go from $u$ to $w$ by a sequence of multiplications by such $s_\delta$,
	the claim follows.

                Now let $\block_1,\ldots,\block_t$ be the blocks of $\mon$ such
                that $w(\block_j)=\beta_j$.  To get a $v$-chain $\alpha_1>\ldots>
                \alpha_t$ in $\mon$ with the desired property,   we start by getting an
                appropriate $\alpha_t$ in $\block_t$.  Let $a$ be the least row index
                of an element of $\block_t$ subject to $D\leq a$.   Choose $(a,b)$
                in $\block_t$ with $b$ being least possible.  We claim that $b<D$.
                If $(a,b)$ is the first element in
	the arrangement of elements of $\block_t$ in non-decreasing order
	of row and column indices,  then $b=c_t$ and so $b<D$.
	If not,  then we may assume that the element $(a',b')$
	immediately to the left of $(a,b)$ is  distinct from
	$(a,b)$, and so by choice of $(a,b)$ we have $a'<D$.
	But  also $b<a'$ because both $(a',b')$ and
	$(a,b)$ are in $\block_t$, and the claim follows.
                Set $\alpha_t=(a,b)$.

                        Now, by Corollary~\ref{cltwothree}, there exists $\alpha_{t-1}$
                        in $\block_{t-1}$ such that $\alpha_{t-1}>\alpha_t$,
                $\alpha_{t-2}$ in $\block_{t-2}$ such that $\alpha_{t-2}>\alpha_{t-1}$,
                and so on.  Writing $\alpha_j=(a_j,b_j)$,   we have
               \begin{eqnarray}
                       b_1<\ldots<b_t<D\leq a_t<\ldots<a_1    \label{epiprooftwo} \end{eqnarray}
               and $a_t=D$ if and only if $D$ is not an entry of $v$.

               Set $y:=s_{\alpha_1}\cdots s_{\alpha_t}v$.     We show $y_j=w_j$
               by showing that $y_j=D=u_j$.     Note first that $D$ occurs as an entry in
               $y$ for similar reasons that $D$ occurs as an entry of $u$.  Now it follows
               from equations (\ref{epiproofone}) and (\ref{epiprooftwo}) that the number
               of entries in $u$ that are less than $D$ equals the corresponding number
               for $y$ (both equal the corresponding number in $v$ reduced by $t$),
               and so $D$ occurs in the same place in $y$ as it does in $u$.

\subsection{The description of $\phi$}\label{ssphi}
We now specify the map $\phi$. Let
$(w,\min)$ be a pair consisting of an element $w$ of $\idn$ with
$w\gneq v$  and a $w$-dominated monomial $\min$,  possibly empty, in $\pos^v$.

Let $\mon_w$ be the monomial in $\pos^v$ associated to $w$ as in Proposition~\ref{psw},
and $k$ be the maximum length of a $v$-chain in $\mon_w$.
For a positive integer $j$, $1\leq j\leq k$,  let $\monwj$
be the subset of elements of $\monw$ that are $j$-deep.
 Letting $w^j$ be the element associated by Proposition~\ref{psw} to $\mon_w^j$,
we have $w=w^1\geq\ldots w^k\geq v$ (see Remark~\ref{rsw}).
Clearly
$\monw=\monwone\supseteq\monwtwo\supseteq\ldots\supseteq\monwk$.

For a positive integer $j$, $1\leq j\leq k$,  let $\minwj$ be the subset of $\min$
of elements $\beta$ such that $\beta$ is the head of a $v$-chain in $\min$ dominated
by $w^j$ but not $w^{j+1}$  (we set $w^{k+1}=v$),  and every $v$-chain in $\min$
with head $\beta$ is dominated by $w^j$.   The $\minwj$ form a partition of
$\min$ (some of the $\minwj$ could be empty).    Two distinct elements belonging
to the same $\minwj$ are not comparable:  if $\beta>\beta'$  and $\beta'$ belongs
to $\minwj$,   choose a $v$-chain $\beta'>\ldots$ in $\min$ that is not dominated by
$w^{j+1}$, so that, by Corollary~\ref{cldom}~(1), the $v$-chain $\beta>\beta'>\ldots$ is not dominated
by $w^j$.

We now further partition each $\minwj$ into subsets called {\em pieces} as follows.
Let $\mon_{w,j}$ denote the set of elements of $\monw$ that are $j$-deep but not
$(j+1)$-deep.    Clearly no two distinct elements of $\mon_{w,j}$ are comparable.
\blemma\label{pphi}
For an element $(r,c)$ of $\minwj$,   there exists a unique element $(R,C)$ of
$\mon_{w,j}$ such that $C\leq c$ and $r\leq R$.
\elemma
\bproof
Let $(r,c)=(r_1,c_1)>\ldots>(r_k,c_k)$ be a $v$-chain in $\min$
that is dominated by $w^j$ but not $w^{j+1}$.    By Lemma~\ref{ldom},  there
exists a $v$-chain $(R_1,C_1)>\ldots>(R_k,C_k)$ of elements of $\monwj$
with $C_i\leq c_i$ and $r_i\leq R_i$.     Note that $(R_1,C_1)$ does not belong to
$\mon_w^{j+1}$, for the former $v$-chain is not dominated by $w^{j+1}$ (we are
using Lemma~\ref{ldom} again here).   This proves
existence.

For uniqueness,  let $(R,C)$ and $(R',C')$ be distinct elements of
$\mon_{w,j}$.   Since $\mon_{w,j}$ is distinguished and no two of
its distinct elements are comparable,   it follows that either
$C<R<C'<R'$ or $C'<R'<C<R$.    But then if $C\leq c<r\leq R$ and
$C'\leq c<r\leq R'$,    there is an obvious contradiction.\eproof

The pieces of $\minwj$ are indexed by $\mon_{w,j}$.    Let  $\beta=(R,C)$ in $\mon_{w,j}$.
The corresponding piece $\piece_\beta$ of $\minwj$ consists,
by definition, of all those $(r,c)$ in $\minwj$
with $r\leq R$ and $C\leq c$.   Of course, some of these pieces could be empty.
Let us arrange the elements of $\piece_\beta$ in non-decreasing order of the row entries;
among those with equal row entries the arrangement is by non-decreasing order of
column entries.     Since no two distinct elements of $\minwj$ are comparable,
the column entries are also now in non-decreasing order.   Suppose the arrangement is
\begin{eqnarray}\label{epiece}
(r_1,c_1), (r_2,c_2),\ldots, (r_p,c_p)		\end{eqnarray}
Note that $C\leq c_1$,  and $r_p\leq R$.     Let $\piece_\beta^*$ denote the monomial
\begin{eqnarray}\label{epiecestar}
	\left\{(r_1,C), (r_2,c_1),\ldots, (r_p,c_{p-1}), (R,c_p)\right\}.
					\end{eqnarray}
 Set
\[
(\min^w_j)^*:=\bigcup\limits_{\mon_{w,j}}\piece_\beta^*
\quad\quad\mbox{ {\rm and} } \quad\quad
\phi(w,\min):=\bigcup\limits_j 
(\min^w_j)^*. 		\]
This finishes the description of $\phi$.

The following lemma
will be used in the proof of Proposition~\ref{ppiphi}.  Its
statement and proof mirror those of Lemma~\ref{lssprime}.
\blemma\label{lppstar}
No two distinct elements of $\min^j_w\cup(\min_w^j)^*$ are comparable.
  In particular, for a given $\beta$ in $\mon_w$, no two distinct
elements of $\piece_\beta\cup\piece_\beta^*$ are comparable.
\elemma
\bproof
We can put the elements in (\ref{epiece}) and (\ref{epiecestar})
together in non-decreasing order of both row and column numbers:
\begin{eqnarray*}
(r_1,C), (r_1,c_1), (r_2,c_1), (r_2,c_2), \ldots,\\
(r_{p},c_{p-1}), (r_p,c_p),  (R,c_p)		\end{eqnarray*}
This means that distinct  elements belonging to $\piece_\beta
\cup\piece_\beta^*$,
where $\piece_\beta$ is a piece of $\min_w^j$,  are not comparable.
If $\alpha=(R,C)$ and $\beta=(S,T)$ are distinct elements
of $\mon_{w,j}$,  then, since $\mon_{w,j}$ is distinguished and no two
of its distinct elements are comparable,
we may assume without loss of generality that
$C<R<S<T$.    If $(a,b)$ belongs to $\piece_\alpha\cup\piece_\alpha^*$ and
$(c,d)$ to $\piece_\beta\cup\piece_\beta^*$,  then
\[	C\leq b< a\leq R<S\leq d<c\leq T	\]
so that $(a,b)$ and $(c,d)$ are not comparable, and we are done.
Note that this shows also that $\piece_\alpha\cup\piece_\alpha^*$
and $\piece_\beta\cup\piece_\beta^*$ are disjoint.\eproof

\subsection{Proof of Proposition~\ref{ppiphi}} \label{ssphipi}
We now turn to the proof of Proposition~\ref{ppiphi}.
We first show that $\phi\circ\pi$ is the identity.

Let $\mon$ be a monomial in $\pos^v$,  and $(w,\mon')$ its
image under $\pi$.    To show that $\phi(w,\mon')=\mon$,
it suffices to prove the following claim:   for $\beta\in\monw$,
if $\block_\beta$ is the block of $\mon$ such that $w(\block_\beta)
=\beta$,  and $\piece_\beta$
the piece of $(w,\mon')$ corresponding to $\beta$,  then
$\block'_\beta=\piece_\beta$.    Given the claim,   we evidently
have $\piece_\beta^*=\block_\beta$, and so
\[
	\phi(w,\mon')=\bigcup\limits_{\beta\in\monw}\piece_\beta^*
		=  \bigcup\limits_{\beta\in\monw}\block_\beta = \mon.	\]

To prove the claim,  it suffices to show that $\block'_\beta\subseteq
\piece_\beta$,  for then, since $\{\piece_\beta\st\beta\in\monw\}$
is a partition of $\mon'=\cup_\beta\block'_\beta$,  equality follows.

Let $(r,c)$ be an element of $\block'_\beta$.  Write $\beta=(R,C)$,
and let $j$ denote the depth of $\beta$ in $\monw$.  It follows from the
definition of $\block'_\beta$ that $C\leq c$ and $r\leq R$, and
so it suffices to show that $(r,c)$ belongs to $\minwj$, where
we write $\min$ instead of $\mon'$ for notational convenience.

\bremark\label{rphipi}
By Corollary~\ref{cltwofour},  $\block_\beta\subseteq\mon_j$ so that
$(r,c)$ belongs to $\mon_j'$.
\eremark

We first show 
that every $v$-chain in $\min$ with head $(r,c)$ is dominated
by $w^j$.    Let $(r,c)=(r_1,c_1)>\ldots>(r_h,c_h)$ be any $v$-chain in $\min$
with head $(r,c)$.    Letting $\beta_i=(R_i,C_i)$ be the element of $\monw$ such that
$(r_i,c_i)$ belongs to $\block_{\beta_i}'$ (see Corollary~\ref{clssprime}), it follows from
the definition of $\block_{\beta_i}'$ that
$C_i\leq c_i$ and $r_i\leq R_i$.
By Corollary~\ref{cltwo}~(2),  we have $\beta=\beta_1>\ldots>\beta_h$.
Since $\beta$
is $j$-deep,  so are $\beta_2,\ldots,\beta_h$,  and thus $\beta_1,\ldots,\beta_h$
belong to $\monwj$.   
By Lemma~\ref{ldom},  we conclude that $w^j$ dominates the $v$-chain
above with head $(r,c)$.

We now construct a $v$-chain with head $(r,c)$ that is not
dominated by $w^{j+1}$.
For this we prove the following lemma.
\blemma
If $\delta=(A,B)$  belongs to $\mon_{w,j+1}$  with  $B\leq c$ and $r\leq A$.
then there exists $(a,b)\in\block_\delta'$ such that $(r,c)>(a,b)$.
\elemma
\bproof
By Corollary~\ref{cltwofour}, $\block_\delta\subseteq\mon_{j+1}$.
Since $w(\block_\delta)=\delta=(A,B)$,  three exists an element of $\mon_{j+1}$
of the form $(A,C)$.    We must have $c<C$,  for otherwise choosing $(x,y)$ in
$\mon_j$ with $(x,y)>(A,C)$,  we have $(x,y)>(r,c)$,  but distinct elements of
$\mon_j\cup\mon_j'$ are not comparable (Lemma~\ref{lssprime}).    Let $(A'b)$ be
the
element of $\mon_{j+1}$ such that $b$ is the least with the property that
$c<b$, and among those with column index $b$ the least row index
possible is $A'$.   In the arrangement of elements of $\mon_{j+1}$ in
non-decreasing order of row and column entries,  there is a portion that
looks like this:
\[
		\ldots,(a,d),(A',b),\ldots		\]
Since there exists an element in $\block_\delta$ (and so in $\mon_{j+1}$)
with column index $B$,   and we have $B\leq c\leq b$,  it follows that $(A',b)$ is
not the left most element, and so we may assume that $(a,d)\neq(A',b)$.  By choice
of $(A',b)$,  we have $d\leq c$.  If now $r\leq a$, choosing
$(x,y)$ in $\mon_j$ with $(x,y)>(a,d)$,   we get $(x,y)>(r,c)$,  but no two distinct
elements of $\monj\cup\mon_j'$ are comparable (Remark~\ref{rphipi} above and
Lemma~\ref{lssprime}).   Thus $a<r$,  so that the element
$(a,b)$ of $\mon_{j+1}'$ satisfies $(r,c)>(a,b)$.\eproof

We now proceed by decreasing induction on $j$ to obtain
a $v$-chain with the desired property.    If $j=k$  (where $k$ is the
depth of $\monw$),  then $w^{j+1}=v$,  and $(r,c)$ is not dominated by $v$.
For $j<k$,  if a $\delta$ as in the hypothesis of
the lemma does not exist,  then $(r,c)$
by itself is not dominated by $w^{j+1}$.   Assuming such a $\delta$ exists,
let $(a,b)$ be as in the conclusion of the lemma.
By induction there exists a $v$-chain $(a,b)>\ldots$ in $\mon'$ that is not
dominated by $w^{j+2}$.   The $v$-chain $(r,c)>(a,b)>\ldots$ is then not
dominated by $w^{j+1}$ (see Corollary~\ref{cldom}~(1)),
and the proof that
$\phi\circ\pi$ is the identity is complete.

We now turn to showing that $\pi\circ\phi$ is the identity.
Let $(w,\min)$ be a pair consisting of an element $w\gneq v$
of $\idn$ and a $w$-dominated monomial $\min$ in $\pos^v$.
Let $\mon=\phi(w,\min)$.     To prove that $\pi(\mon)=(w,\min)$,
it suffices to show that, for any piece $\piece$ of
$(w,\min)$, $\piece^*$ is a block of $\mon$.
This in turn follows from the following two assertions (recall that, for
$\beta\in\monw$, $\piece_\beta$ denotes the piece of $(w,\min)$
indexed by $\beta$):
\begin{enumerate}
	\item
		For $\beta\in\monw$,  there exists a block
		$\block$ of $\mon$ such that
		$\piece_\beta^*\subseteq\block$ (this block is
		then uniquely determined, for $\piece_\beta^*$
		is non-empty and the blocks are disjoint).
	\item
		If $\beta\neq\beta'$ are elements of $\monw$,
                                and $\block$, $\block'$  are the blocks
		determined as in (1) respectively by $\beta$
		and $\beta'$, then $\beta\neq\beta'$. 
\end{enumerate}
To see why it suffices to prove the above assertions, we fix an element
$\beta$ of $\mon_w$,  let $\piece_\beta$ be the piece of $(w,\min)$
corresponding to $\beta$,   let $\block$ be the block of $\mon$ as in (1),
and show that $\piece_\beta=\block$.
 Given $(r,c)\in\block$,  since $\mon:=\cup_{\beta\in\monw}
\piece_\beta^*$ by definition of $\phi$,  there is $\beta'$ in
$\monw$ such that $(r,c)\in\piece_{\beta'}^*$, and by (1) there
is a $\block'$ such that $\piece_{\beta'}^*\subseteq\block'$,
but then $\block=\block'$ since $\block\cap\block'$ contains
$(r,c)$ and is therefore non-empty (the blocks form a partition
of $\mon$),  which means $\beta=\beta'$ by (2), so that
$(r,c)\in\piece_\beta^*$.

To prove the assertions above,  we make use of the following
lemma:
\blemma\label{lfour}
Let $\beta$ and $\beta'$  be elements of $\monw$, and $(r,c)$
an element of $\piece_\beta^*$.   Then there exists $(r',c')$
in $\piece_{\beta'}^*$  with $(r',c')>(r,c)$ if and only if $\beta'>\beta$.
\elemma
\bproof
Let us write $\beta=(R,C)$ and $\beta'=(R',C')$.
Let $j$ and $j'$ denote respectively the depths of $\beta$
and $\beta'$ in $\monw$.
It is clear from the definition of $\piece_\beta^*$ that $r\leq R$ and $C\leq c$.

First suppose that $\beta'>\beta$. We then have
$R<R'$, $C'<C$,   and $r\leq  R$, $C\leq c$ , 
so that $r<R'$.
 There exists an element in $\piece_{\beta'}^*$ with row
index $R'$,  say $(R',a)$.   Let $(r',c')$ be the  element
of $\piece_{\beta'}^*$ such that $r'$ is least with the property
that $r<r'$,   and among those with row index $r'$ the least possible
column index is $c'$.  If $c'<c$,  then $(r',c')$ has the desired
property.   So it is enough to assume that $c\leq c'$
and arrive at a contradiction.

In the arrangement of elements of $\piece_{\beta'}^*$ in non-decreasing
order of row and column indices,  there is a portion that looks like
this:
\[	\ldots,(a,b),(r',c'),\ldots		\]
Since $C'<c'$  (for $C'<C\leq c\leq c'$) and the first element of
$\piece_{\beta'}^*$ has $C'$ for its column index,  it follows
that $(a,b)$ exists,  and so we may assume that $(a,b)\neq(r',c')$.
By choice of $(r',c')$,  it follows that $a<r'$, and further
that $a\leq r$.

We now set things up for an application of Corollary~\ref{cldom}.  To this
end, note that $(a,c')$ is an element of $\piece_{\beta'}$ with
$C\leq c\leq c'$ and $a\leq r\leq R$.   Since $\beta'>\beta$,  we have $j'<j$.  Since
$(a,c')$ belongs to $\piece_{\beta'}\subseteq\min_{j'}^w$,  there
exists a $v$-chain
$(a,c')=\alpha_1>\alpha_2>\ldots$ in $\min$ with head $(a,c')$
that is dominated by $w^{j'}$ but not by $w^{j'+1}$.
Now, by Corollary~\ref{cldom},  it follows that
$\alpha_2>\alpha_3>\ldots$
is dominated by $w^{j'+1}$ but not by $w^{j'+2}$,
$\alpha_3>\alpha_4>\ldots$
is dominated by $w^{j'+2}$ but not by $w^{j'+3}$, and so on,
so that
$\alpha_{j-j'+1}>\ldots$
is dominated by $w^{j}$ but not by $w^{j+1}$.
Further,  every $v$-chain  with head $\alpha_{j-j'+1}$ is dominated
by $w^j$, for otherwise by prefixing $\alpha_1>\ldots>\alpha_{j-j'}$ we would
get a $v$-chain with head $(a,c')$ that is not dominated by $w^{j'}$ (again using
Corollary~\ref{cldom}), a contradiction.
Writing $\alpha_{j-j'+1}=(x,y)$,  we have $C\leq c\leq c'<y$
and $x<a\leq r\leq R$,  so that $\beta>(x,y)$.  Thus $(x,y)$ belongs
to $\piece_\beta$.  We also have $(r,c)>(x,y)$,  but this is a
contradiction since distinct elements of $\piece_\beta\cup
\piece_{\beta}^*$ are not comparable (Lemma~\ref{lppstar}),
and the proof of the ``if'' part is complete.

For the converse,  let $(r',c')$ be in $\piece_{\beta'}^*$
with $(r',c')>(r,c)$.   
Since $\beta$ and $\beta'$ belong to $\mon_w$ and $\mon_w$
is distinguished, it follows
that one of the following four possibilities holds:
\begin{eqnarray*}
	 C'<C<R<R', &  \\
	C<C'<R'<R,  \quad
	C<R<C'<R', & \quad
	C'<R'<C<R,  		\end{eqnarray*}
What we want to show is precisely that the first possibility holds.
It therefore suffices to rule out the last three possibilities.
In case (3),   we have $c<r\leq R<C'$ and $C'\leq c'<c$,
a contradiction.   In case (4),  we have $c<r< r'\leq
R'<C$ and $C\leq c$,  a contradiction.

Now suppose that $C<C'<R'<R$.  We have $j'>j$ since
$\beta>\beta'$.  Since
$r<r'\leq R'<R$,  it follows from the way $\piece_\beta^*$ is
obtained from $\piece_\beta$ that there exists an element
$(r,a)$ in $\piece_\beta$ with
$c\leq a$.
Let $(r,a)=\alpha_1>\alpha_2>\ldots$ be a $v$-chain in $\min$
with head $(r,a)$
that is dominated by $w^{j}$ but not by $w^{j+1}$.
Note that $C'\leq c'<c\leq a$ and $r<r'\leq R'$.
By Corollary~\ref{cldom}, it follows that
$\alpha_2>\alpha_3>\ldots$
is dominated by $w^{j+1}$ but not by $w^{j+2}$,
$\alpha_3>\alpha_4>\ldots$
is dominated by $w^{j+2}$ but not by $w^{j+3}$, and so on,
so that
$\alpha_{j'-j+1}>\ldots$
is dominated by $w^{j'}$ but not by $w^{j'+1}$.
Now $\alpha_{j'-j+1}$ belongs to $\piece_{\beta'}$,
and $(r',c')>(r,a)>\alpha_{j'-j+1}$,  but this is a
contradiction since no two distinct elements of $\piece_{\beta'}\cup
\piece_{\beta'}^*$ are comparable (Lemma~\ref{lppstar}).\eproof

It follows immediately from the lemma that $\piece_\beta^*
\subseteq\monj$ where $j$ is the depth of $\beta$ in $\monw$.
That all of $\piece_\beta^*$ is contained in a single block
of $\monj$ follows immediately from the definition of
$\piece_\beta^*$.   This proves (1).

Let us now prove (2).   Suppose that $j$ and $j'$ are
respectively the depths of $\beta$ and $\beta'$ in $\monw$.
From the lemma we have $\piece_\beta^*\subseteq\monj$ and
$\piece_{\beta'}^*\subseteq\mon_{j'}$, so $\block\subseteq\monj$
and $\block'\subseteq\mon_{j'}$.
Suppose that $j\neq j'$.
Then, since
$\block$ and $\block'$ are non-empty, and $\monj$ and
$\mon_{j'}$ are disjoint, it follows that
$\block\neq\block'$.    So suppose that $j=j'$.
Writing $\beta=(R,C)$ and $\beta'=(R',C')$,  we may
assume without loss of generality that $C<R<C'<R'$.
For $(r,c)$ in $\piece_\beta^*$ and $(r',c')$ in $\piece_{\beta'}^*$,
we have $r\leq  R<C'\leq c'$, so that,
by the definition of the partition of $\mon_j$ into blocks, $\block\neq\block'$.
This finishes the proof that $\phi\circ\pi$ is the identity
and so also that of Proposition~\ref{ppiphi}.

\section{Interpretations}\label{sinterpretation}
Fix elements $v,w$ in $\idn$ with $v\leq w$. 
It follows from Corollary~\ref{cmain} that the
multiplicity of the Schubert variety $X_w$ in
the Grassmannian $\gdn$ at the point $e^v$ can be
interpreted as the cardinality of a certain set
of non-intersecting lattice paths.    We first
illustrate this by means of two examples and
then provide a justification for the interpretation.
In the case $v=(1,\ldots,d)$,  this interpretation
is due to Herzog and Trung \cite[Theorem~3.3]{ht}.

\bexample\label{eone}
Let $d=14$, $n\geq 27$,  
\begin{eqnarray*}
v&=&(1,2,3,4,5,10,11,12,13,18,19,20,21,22),
\mbox{\begin{rm} and \end{rm}}\\
w&=&(1,4,5,9,12,13,16,17,19,22,24,25,26,27)\end{eqnarray*}
so that 
\[
\mon_w=\left\{(9,3), (16,11), (17,10), (24,21), (25,20),  (26,18),
			(27,2) \right\}	\]

\begin{figure}\label{feone}
\begin{center}
\mbox{\epsfig{file=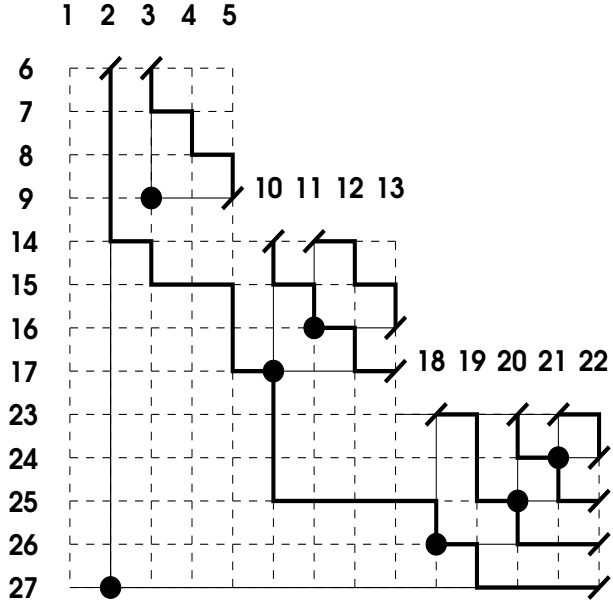,height=8cm,width=8cm}}
\caption{A tuple of non-intersecting lattice paths as in 
	Example~\ref{eone}}
\end{center}
\end{figure}

The grid depicting the points of $\pos^v$ is shown
in Figure 4.  The solid dots represent the
points of $\mon_w$.   From each point $\beta$
of $\mon_w$ we
draw a vertical line and a horizontal line.
Let $\bstart$ and $\bfinish$ denote respectively the points
where the vertical line and the horizontal line meet
the boundary.     For example, $\bstart=(14,11)$ and
$\bfinish=(16,13)$ for $\beta=(16,11)$.

A lattice path between a pair of such points $\bstart$
and $\bfinish$ is a sequence $\alpha_1,\ldots,\alpha_q$
of elements of $\pos^v$ with $\alpha_1=\bstart$,
$\alpha_q=\bfinish$,  and for $j$, $1\leq j\leq q-1$,
writing $\alpha_j=(r,c)$,  $\alpha_{j+1}$ is either
$(r+1,c)$ or $(r,c+1)$.     If  $\bstart=(r,c)$ and
$\bfinish=(R,C)$,   then $q=(R-r)+(C-c)+1$.

Write $\mon_w=\{\beta_1,\ldots,\beta_p\}$.  Now consider
the set of all $p$-tuples of paths $(\path_1,\ldots,\path_p)$,
where $\path_j$ is a lattice path between $\bstartj$
and $\bfinishj$,  and no two $\path_j$ intersect.
A particular such $p$-tuple is shown in Figure 4.
The number of such $p$-tuples is the multiplicity of 
$X_w$ at the point $e^v$.\eexample

\bexample\label{etwo}
Let us draw, in a simpler case,
the pictures of all possible tuples of non-intersecting
lattice paths as defined in the above example. 
Let $d=6$, $n\geq 13$, $v=(1,2,3,8,9,10)$, and
$w=(4,6,7,10,11,13)$,  so that $\mon_w=\{(4,3), (6,2),
(7,1), (11,9), (13,8)\}$.   Figure 5 shows all
the $5$-tuples of 
non-intersecting lattice paths.  There are $9$ of them
and thus the multiplicity in this case is $9$.\eexample

\begin{figure}\label{fetwo}
\begin{center}
\mbox{\epsfig{file=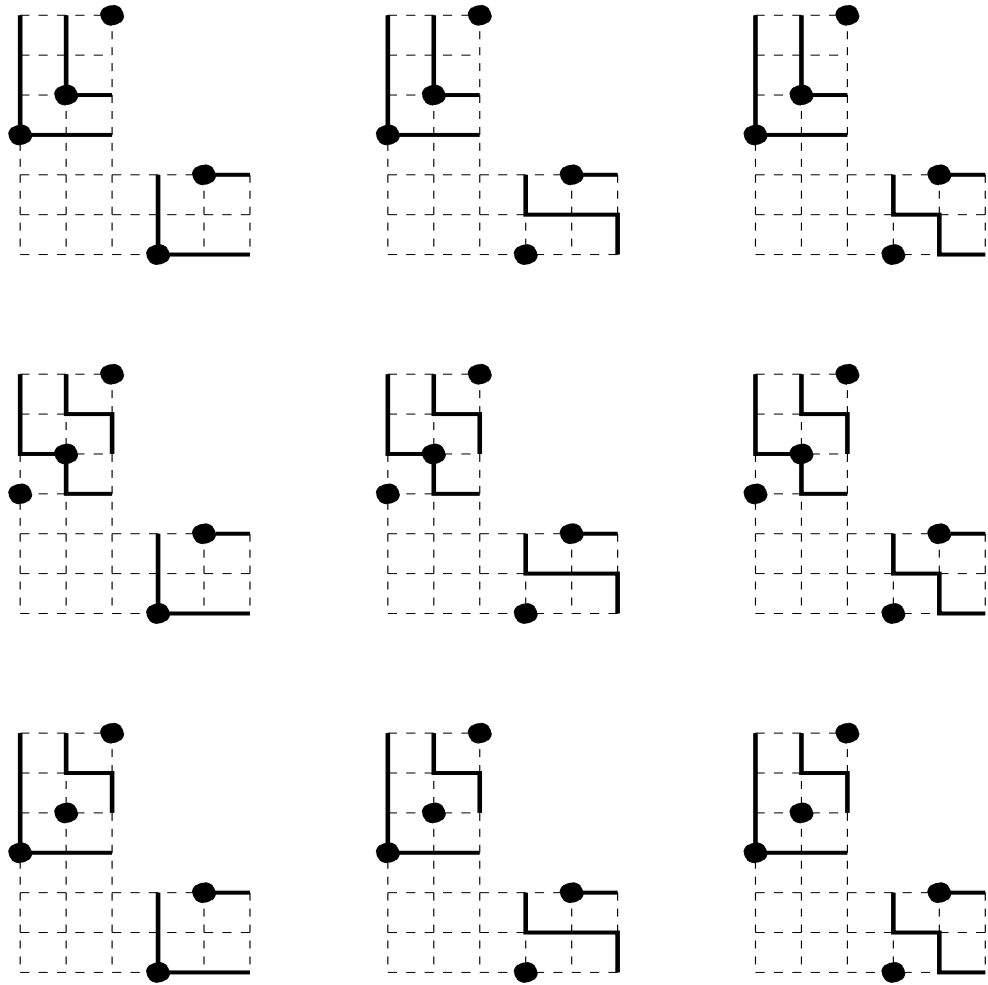,height=13cm,width=13cm}}
\caption{All the tuples of non-intersecting lattice paths as in 
	Example~\ref{etwo}}
\end{center}
\end{figure}

Let us now see why this interpretation is justified.
Denote by $\uvw$ the set of tuples of non-intersecting lattice
paths corresponding to the pair $(v,w)$,   by $\nvw$ the set of
maximal square-free $w$-dominated monomials in $\pos^v$.
There is an obvious natural map from $\uvw$ to monomials in $\pos^v$.
In fact,  each element of $\uvw$ can be thought of as a monomial:
each lattice path is a subset of $\pos^v$ and as such is a 
monomial in $\pos^v$;     we identify an element $\path
=(\path_1,\ldots,\path_p)$ of $\uvw$ with the monomial
$\path_1\cup\cdots\cup\path_p$. 
We denote the monomial also by $\path$.
Since the paths are non-intersecting,
$\path$ is square-free.   We will argue below that $\path$
belongs to $\nvw$ and that this association is a bijection.
The idea is that
the block decomposition of an element of $\nvw$ identifies it
as an element of $\uvw$.
Since all monomials $\path$ have clearly the same cardinality,
the same will be true for elements of $\nvw$.

Let us for the moment take for granted this bijection between
$\uvw$ and $\nvw$ and justify the interpretation.
Denote by $\rvw$ the set of maximal square-free $w$-dominated
monomials in $\roots^v$.
%
If $\mon$ is an element of $\rvw$ and $(r,c)$ an element of
$\roots^v\setminus\pos^v$ not belonging to $\mon$,  then, adding
$(r,c)$ to $\mon$ would still keep it square-free and $w$-dominated,
violating its maximality.  Thus $\mon\supseteq\roots^v\setminus
\pos^v$,  and the map $\mon\mapsto\mon\cap\pos^v$ is a bijection
between $\rvw$ and $\nvw$.   Combining this with the bijection
$\uvw\cong\nvw$,  we get a natural bijection $\uvw\cong\rvw$.
Furthermore, the elements
of $\rvw$ all have the same cardinality.  The justification
is now completed by an appeal to Corollary~\ref{cmain}.

We now sketch a proof that $\uvw$ and $\nvw$ are naturally bijective.
We first show that $\path$ is $w$-dominated.   It is not hard to see
that $\path=(\path_1,\ldots,\path_p)$ is the block decomposition
of $\path$.      Thus $w$ is the first component of $\pi(\path)$,
where $\pi$ is the map described in subsection~\ref{sspi}.   By
Proposition~\ref{ppi}~(4),   $w$  dominates $\path$.

To show the maximality of $\path$ with respect to being square-free
and $w$-dominated,  let $(r,c)$ be an element of $\pos^v\setminus\path$.
Consider the block decomposition of $\path^{\#}:=\path\cup\{(r,c)\}$.
There is no way in which $(r,c)$ can fit into some block 
$\path_j$ without affecting $w(\path_j)$.     Thus the first component
of $\pi(\path^{\#})$ is strictly larger than $w$,  and  it  follows
from Proposition~\ref{ppi}~(4) that $\path^{\#}$ is not $w$-dominated.

We have thus shown that elements of $\uvw$ when thought of
as monomials belong to $\nvw$.   We will now show that each
element of $\nvw$ arises this way.

Let $\mon$ be an element of $\nvw$ and $\block$ a block of $\mon$.
Suppose that $(r,c)$ and $(R,C)$ are consecutive members in
the arrangement of elements of $\block$ by non-decreasing
row and column indices.    Then either $r<R$ or $c<C$ since
$\mon$ is square-free.

Suppose that $r<R$.   If $c<C$,  then consider $\mon\cup\{(r,c)\}$.
It is not hard to see that $\mon\cup\{(r,c)\}$ has the same
block decomposition as $\mon$ except that $\block$ is replaced
by $\block\cup\{(r,c)\}$.  This violates maximality of $\mon$.
Thus $c=C$.

If $r<R-1$, consider $\mon\cup\{(R-1,C)\}$.  This again has the same
block decomposition as $\mon$ except that $\block$ is replaces by
$\block\cup\{(R-1,C)\}$.   This violates the maximality of $\mon$.
Thus $r=R-1$.

In a similar fashion,   we can argue that if $c<C$,  then $r=R$ and
$c=C-1$.    Furthermore,  we can show similarly that
the first member of $\block$ is $(a,b)$ and the last $(A,B)$,
where $w(\block)=(A,b)$,  $a$ is the least element of $\{1,\ldots,n\}
\setminus\{v_1,\ldots,v_d\}$ that is bigger than $b$,  and $B$ is
the biggest element of $\{v_1\ldots,v_d\}$ that is smaller than $A$.

The block decomposition of $\mon$ thus gives a tuple of non-intersecting
lattice paths that belong to $\uvw$.

\subsection{The Gr\"obner basis interpretation}\label{ssgroebner}
Fix elements $v,w$ in $\idn$ with $v\leq w$.    As in the case
$v=(1,\ldots,d)$ considered in \cite{kl},   the main theorem has
an interpretation in terms of Gr\"obner basis.

Recall from \S\ref{ssmt} that the ideal $\left(f_\theta=p_\theta/p_v\st
\mbox{\begin{rm} $\theta$ in $\idn$, $\theta\nleq w$\end{rm}}\right)$
in the polynomial ring $P:=k[X_\beta\st\beta\in\roots^v]$ defines the
tangent cone to the Schubert variety $X_w$ at the point $e^v$.
The $f_\theta$ are homogeneous polynomials in the variables
$X_\beta$, $\beta\in\roots^v$.        The following lemma shows that
a subset of these generators suffices:
\blem\label{lgenerator}
The ideals $(f_\theta\st\theta\in\idn,\theta\nleq w)$ and $(f_\theta\st
\theta\in\idn, v\leq\theta\nleq w)$ are equal.
\elem
\bproof
Let $\theta$ be an element of $\idn$ with $\theta\nleq w$.   We will show
that $f_\theta$ belongs to $(f_\mu\st\mu\in\idn, v\leq\mu\nleq w)$.
Fix $k$, $1\leq k\leq d$, such that $\theta_k>w_k$.     Let $a$ be the largest
integer such that $1\leq a\leq d$ and $v_a\leq w_k$.    Consider the
elements $\mu$ of $\idn$ of the form
\[
        \left\{v_{i_1},\ldots,v_{i_{k-1}}\right\}\cup
      \left\{v_{a+1},\ldots,v_d\right\}\cup
      \left\{\theta_{j_1},\ldots,\theta_{j_{a-k+1}}\right\}\]
where $1\leq i_1<\ldots<i_{k-1}\leq a$, $k\leq j_1<\ldots<j_{a-k+1}\leq n$,
and
 $ \left\{v_{a+1},\ldots,v_d\right\}$ is disjoint from
 $\left\{\theta_{j_1},\ldots,\theta_{j_{a-k+1}}\right\}$.

 Let us first observe that such a $\mu$ satisfies $v\leq \mu\nleq w$.
 For $1\leq j\leq k-1$,   we have $\mu_j=v_{i_j}\geq v_j$  (since $i_j\geq j$).
 For $k\leq j\leq a$,  we have $\mu_j\geq \mbox{\begin{rm}min\end{rm}}
 \left\{v_{a+1},\theta_k\right\}>w_k
 \geq v_a$.   (This also shows that $\mu_k>w_k$.)   Let now $a<j$.   Consider
 the subset of the first $j-k+1$ elements of $\{v_{a+1},\ldots,v_d\}\cup
 \{\theta_{j_1},\ldots,\theta_{j_{a-k+1}}\}$ when the elements are arranged in
 increasing order.    This subset contains $v_j$, and since $\mu_j$ is the largest
 element of this subset,   $\mu_j\geq v_j$.    This finishes the proof that $\mu$
 satisfies $v\leq \mu\nleq w$.

 Let $M$ denote the sub-matrix of a matrix as in Figure~1 consisting of
 the rows $\theta_1,\ldots,\theta_d$.    Let $N$ be the $(d-k+1)\times a$
 sub-matrix of $M$ determined by the rows $\theta_k,\ldots,\theta_d$ and
 the first $a$ columns.    Each $f_\mu$ is the determinant of an $(a-k+1)\times
 (a-k+1)$ sub-matrix of $N$.   And each such non-zero determinant of
 a sub-matrix of $N$ equals $f_\mu$ for some $\mu$.

 That $f_\theta$ belongs to the ideal $(f_\mu\st \mu \mbox{\begin{rm}
  as above\end{rm}})$
 follows from
 the following fact about matrices,  which is easily seen (by the Laplace expansion,
 for example):  Let $d,a,k$ be integers such that $1\leq k\leq a\leq d$.
 Let $M$ be a generic matrix of indeterminates of size $d\times d$.     Let $N$
 be a sub-matrix of $M$ of size $(d-k+1)\times a$  (in other words, $N$ is obtained
 by specifying $(d-k+1)$ row indices and $a$ column indices).     Then, in
 the polynomial ring obtained by adjoining to a field the entries of $M$,
 the determinant of $M$ belongs to the ideal generated by the determinants
 of sub-matrices of $N$ of size $(a-k+1)\times (a-k+1)$.\eproof

 For $\theta$ in $\idn$ with $v\leq\theta$, let $\mon_\theta$ denote the
monomial in $\pos^v$ (and so also in $\roots^v$)
associated to
$\theta$ as in Proposition~\ref{psw}.   We abuse notation and write
$\mon_\theta$ also for the corresponding monomial in the variables
$X_\beta$, $\beta\in\roots^v$.

\bprop\label{pgroebner}
Let $v,w$ be elements of $\idn$ with $v\leq w$.  Fix any term
order on the monomials in the polynomial ring
$k[X_\beta\st\beta\in\roots^v]$ such that
the initial term of $f_\theta$ with respect to this order is
$\mon_\theta$
for $\theta$ in $\idn$ with $\theta\geq v$.    Then
$\left\{f_\theta\st\theta\in\idn, v\leq\theta\nleq w\right\}$ is
a Gr\"obner basis with respect to this order.
\eprop
\bproof
Denoting by $I$ the ideal $\left(f_\theta
\st\theta\in\idn, v\leq\theta\nleq w\right)$,  by $\init{f}$ the
initial term of a polynomial $f$ in this term order,  and by $\init{I}$
the ideal $\left(\init{f}\st f\in I\right)$,  we clearly have a graded
surjection
\[
P/\left(\init{f_\theta}\st\theta\in\idn, v\leq\theta\nleq w\right)
\surjection P/\init{I}
\]
The assertion of the proposition is that this map is an isomorphism.
To prove this,  it is enough to show that both graded rings have the
same Hilbert function.

The Hilbert function of $P/\init{I}$ is the same as that of the tangent cone
$P/I$,  and so by Theorem~\ref{tmain} it equals the cardinality of
$S^v_w(m)$ in degree $m$, where $\svwm$ is the set of $w$-dominated
monomials in $\roots^v$ of degree $m$.  We now show that a monomial
in $\roots^v$ belongs to the ideal $\left(\init{f_\theta}\st\theta\in\idn,
v\leq\theta\nleq w\right)$ if and only if it is not $w$-dominated, and
this will complete the proof.

Suppose that a monomial $\mu$ in $\roots^v$ is not $w$-dominated.
Let $\beta_1>\ldots>\beta_t$ be a $v$-chain of elements in $\mu$
such that $w\ngeq s_{\beta_1}\cdots s_{\beta_t}v$.  The monomial
$\{\beta_1,\ldots,\beta_t\}$ is distinguished and so by
Proposition~\ref{psw} it corresponds to an element $\theta$ in
$\idn$ with $v\leq \theta$.   We have $\theta=s_{\beta_1}
\cdots s_{\beta_t}v\nleq w$ so that $\init{f_\theta}$ divides $\mu$.

For the converse, it is enough to show that $\init{f_\theta}$ is
not $w$-dominated for $\theta\nleq w$, and this follows from the
following lemma.\eproof
\blemma
For $\theta$ in $\idn$ with $v\leq\theta$,  we have $\theta\leq w$
if and only if $w$ dominates $\mon_\theta$.
\elemma
\bproof
Let us consider the image of the monomial $\mon_\theta$ under the
map $\pi$ of subsection~\ref{sspi}.    The blocks of $\mon_\theta$
are precisely its singleton subsets.   Thus $\pi(\mon_\theta)
=(\theta,\mbox{{\rm empty monomial}})$.   Proposition~\ref{ppi}~(4)
now says that $\theta$ is the least element of $\idn$ that dominates
$\mon_\theta$.   The lemma follows.\eproof

We finish by listing some term orders that satisfy the requirement
of Proposition~\ref{pgroebner}.    Fix notation and
terminology as in \S15.2 of \cite{ebud}.  Defined below are
four partial orders $\poone$, $\potwo$, $\pothree$, and $\pofour$
on the elements of $\roots^v$.  We use $>_j$ to denote also any
total order that refines the partial order $>_j$.  The homogeneous
lexicographic order induced by $>_1$ or by $>_2$ and the
reverse lexicographic order induced by $\pothree$ or by $\pofour$
satisfy the requirement of Proposition~\ref{pgroebner}.

The partial orders $\poone$, $\potwo$, $\pothree$, and $\pofour$
on elements of $\roots^v$ are as follows:
\begin{itemize}
\item
	for $(r,c)$ in $\pos^v$ and $(a,b)$ in $\roots^v\setminus\pos^v$,
	$(r,c)>_j(a,b)$ for $j=1,2,3,4$.
\item
	for $(r,c)$ and $(r',c')$ in $\pos^v$,
	\begin{enumerate}
	\item	$(r,c)\poone(r',c')$ if either (a)~$r<r'$ or
			(b)~$r=r'$ and $c>c'$.
	\item	$(r,c)\potwo(r',c')$ if either (a)~$c>c'$ or
			(b)~$c=c'$ and $r<r'$.
	\item	$(r,c)\pothree(r',c')$ if either (a)~$r<r'$ or
			(b)~$r=r'$ and $c<c'$.
	\item	$(r,c)\pofour(r',c')$ if either (a)~$c>c'$ or
			(b)~$c=c'$ and $r>r'$.
			\end{enumerate}
			\end{itemize}

\end{document}